\documentclass[12pt,reqno,a4paper]{amsart}
\usepackage{
    amsmath,  amsfonts, amssymb,  amsthm,   amscd,
    gensymb,  graphicx, comment,  etoolbox, url,
    booktabs, stackrel, mathtools,enumitem, mathdots,  microtype, lmodern,    mathrsfs, graphicx, tikz,  longtable,tabularx, float, tikz, pst-node, tikz-cd, multirow, tabularx, amscd,  bm, array, makecell, diagbox, booktabs,ragged2e, caption, subcaption }
\usepackage{makecell,slashbox}
\usepackage{xcolor}
\usepackage[all]{xy}
\usepackage{graphicx}

\usepackage[utf8]{inputenc}
\usepackage{microtype, fullpage, wrapfig,textcomp,mathrsfs,csquotes,fbb}
\usepackage[colorlinks=true, linkcolor=blue, citecolor=blue, urlcolor=blue, breaklinks=true]{hyperref}
\usepackage[capitalise]{cleveref}
\usepackage{todonotes}
\usetikzlibrary{positioning}
\usetikzlibrary{shapes,arrows.meta,calc}
\usetikzlibrary{arrows}

\newtheorem{theorem}{Theorem}[section]

\newtheorem{corollary}[theorem] {Corollary}
\newtheorem{definition}[theorem]{Definition}
\newtheorem{example}[theorem]{Example}
\newtheorem{lemma}[theorem]{Lemma}

\newtheorem{proposition}[theorem]{Proposition}
\newtheorem{remark}[theorem]{Remark}

\newcommand{\TC}{\mathrm{TC}}

\newcolumntype{x}[1]{>{\centering\arraybackslash}p{#1}}

\begin{document}

\title{Bivariate Topological Complexity:\\ A Framework for Coordinated Motion Planning}

\author{ J.M. Garc\'ia-Calcines and J.A. Vilches   }

\maketitle

\begin{abstract}
We introduce a bivariate version of topological complexity, $\mathrm{TC}(f,g)$, associated with two continuous maps $f\colon X\to Z$ and $g\colon Y\to Z$. This invariant measures the minimal number of continuous motion planning rules required to coordinate trajectories in $X$ and $Y$ through a shared target space $Z$. It recovers Farber’s classical topological complexity when $f=g=\mathrm{id}_X$ and Pavešić’s map-based invariant when one of the maps is the identity.

We develop a structural theory for $\mathrm{TC}(f,g)$, including symmetry, product inequalities, stability properties, and a collaboration principle showing that, when one of the maps is a fibration, the complexity of synchronization is controlled by the other. We also introduce a homotopy-invariant bivariate complexity $\mathrm{TC}_H(f,g)$ of Scott type, defined via homotopic distance, and study its relationship with the strict invariant.

Concrete examples reveal rigidity phenomena with no analogue in the classical case, including strict gaps between $\mathrm{TC}(f,g)$ and $\mathrm{TC}_H(f,g)$ and situations where synchronization becomes impossible. Cohomological estimates provide computable obstructions in both the strict and homotopy-invariant settings. 
%{\red Links with Scott's and Murillo-Wu's are also established, which justifies how bivariate $TC$ generalizes different approaches.}
\end{abstract}

%\vspace{0.5cm}
%\noindent{2010 \textit{Mathematics Subject Classification} : 55M30, 55M15, 55P99.}\\
%\noindent{\textit{Keywords} : {blablabla}}
%\vspace{0.2cm}

%\section*{Introduction}

\section*{Introduction}

The notion of topological complexity, introduced by Farber in the early 2000s  \cite{Far,Far2}, was designed to quantify the minimal number of continuous rules required to perform motion planning in a configuration space. For a topological space $X$, the invariant $\mathrm{TC}(X)$ measures how intricate the task of algorithmically assigning paths between arbitrary pairs of points in $X$ can be. Since its introduction, this invariant has become central in applied and pure topology, with connections ranging from robotics to homotopy theory.  

Later developments extended the scope of topological complexity. In particular, Pave\v{s}i\'c defined in \cite{Pav2} the topological complexity of a map $f \colon X \to Z$, thereby capturing situations in which motion planning depends not only on the internal structure of $X$ but also on an external task space $Z$. Scott subsequently introduced in \cite{Sc} another generalization, denoted $\mathrm{TC}(f)$. Although Scott did not frame his definition in terms of homotopic distance, it can be reinterpreted within this context thanks to the notion introduced later by Mac\'ias--Virg\'os and Mosquera--Lois in \cite{M-M}. This reinterpretation clarifies the position of Scott’s invariant as part of a broader framework and will be developed in the sequel. Related coincidence-type invariants were later considered by Murillo and Wu \cite{M-W}, arising from a different conceptual motivation but leading, in the univariate setting, to numerical invariants that agree with Scott’s.

In many practical and theoretical scenarios, however, multiple subsystems interact through a shared codomain. Think of robotic agents, distributed sensors, or parallel computational processes: each subsystem has its own configuration space but must coordinate through a common target space. Such settings motivate the introduction of a \emph{bivariate topological complexity} associated with two continuous maps $f \colon X \to Z$ and $g \colon Y \to Z$, denoted by $\mathrm{TC}(f,g)$. Given an input pair $(x,y)\in X\times Y$, a local $(f,g)$--motion--planning rule produces two paths $\alpha\in X^I$ and $\beta\in Y^I$ with $\alpha(0)=x$ and $\beta(1)=y$ whose images are synchronized in $Z$, in the sense that $f\circ\alpha=g\circ\beta$. No compatibility condition on $(x,y)$ is assumed a priori; the invariant counts the minimal number of such continuous local rules needed to cover $X\times Y$ (and is set to infinity if such a finite cover does not exist). Conceptually, this extends the motion--planning problem from single systems to the coordination of two independent agents constrained by a common task space, with synchronization enforced along the entire paths rather than only at endpoints. 

From a theoretical perspective, this construction unifies and generalizes existing invariants. When $f=g=\mathrm{id}_X$, one recovers Farber’s original invariant $\mathrm{TC}(X)$. When one of the maps is the identity on $Z$, the invariant specializes to the map-based framework developed by Pave\v{s}i\'c. The bivariate version thus sits at the intersection of these approaches, while opening new territory: it is sensitive to the interplay between different configuration spaces, and it admits both strict and homotopy-invariant formulations. In this sense, $\mathrm{TC}(f,g)$ can be regarded as a numerical measure of how costly it is to synchronize two systems through a common target space, extending the classical homotopical framework of sectional category and pullbacks.  

While the bivariate invariant $\mathrm{TC}(f,g)$ can be motivated as a measure of synchronization cost, it is important to emphasize that what it formally quantifies is the minimal number of continuous local rules required to coordinate trajectories in $X$ and $Y$ whose images under $f$ and $g$ coincide in $Z$. In particular, it measures the topological complexity of the space of such synchronized path-lifts rather than a physical cost in any specific model. For instance, when $g=\mathrm{id}_Z$, the invariant specializes to the map-based framework of Pave\v{s}i\'c, showing that one subsystem ($Y=Z$) simply reproduces the trajectory of the other. Thus the bivariate framework captures the topological constraints of coordination, which may in turn reflect synchronization challenges in applied contexts.  

The range of potential applications is wide. For instance, they include distributed computing, where different processes must agree on a common output, or multi-component mechanical systems, such as exoskeletons or drone swarms, where consistent coordination across subsystems is required.  

Beyond these general areas, two illustrative scenarios may help convey the role of the bivariate invariant:  

\begin{itemize}
  \item \textbf{Robotics (drone--arm coordination).} A drone ($X$) is tasked with transporting a piece while a robotic arm ($Y$) holds a support in a shared workspace ($Z$). The maps $f,g$ encode, respectively, the placement of the piece and the position of the support. A synchronized motion plan requires trajectories that remain compatible in $Z$, reflecting the challenge of coordinating two devices whose actions must fit together continuously.  

  \item \textbf{Sensor fusion (camera--LiDAR).} A camera ($Y$) and a LiDAR ($X$) produce data in a common feature space ($Z$). If their fields of view do not overlap, synchronization becomes impossible, highlighting the rigidity of the strict invariant. In contrast, if they do overlap but one sensor is subject to distortions, the strict version may still fail while the homotopy-invariant version $\mathrm{TC}_H(f,g)$ remains finite. This illustrates the difference between exact calibration (strict) and error-tolerant fusion (homotopical).  
\end{itemize}

These examples illustrate how the bivariate framework can model cooperation and synchronization in realistic settings, while leaving open the precise mathematical mechanisms by which such effects occur. They provide a concrete narrative counterpart to the abstract properties that will be developed in the sequel.  

In this article we introduce the invariant $\mathrm{TC}(f,g)$, establish its first structural properties, and compare it with its homotopy-invariant counterpart $\mathrm{TC}_H(f,g)$. Special emphasis is placed on situations where one of the maps is a fibration, as well as on cohomological lower bounds derived from the structure of $H^*(X)$, $H^*(Y)$, and $H^*(Z)$. Along the way we highlight examples illustrating novel phenomena that have no analogue in the classical case, thereby underscoring the richness of the bivariate framework.

The results of this paper are organized into two closely related but conceptually distinct layers: a strict bivariate theory, measuring exact synchronization through a common target space, and a homotopy–invariant theory of Scott type, which provides a flexible and computable relaxation of the strict framework.

On the strict side, we first introduce the invariant $\mathrm{TC}(f,g)$ and analyze its basic geometric features. In particular, we study the notion of synchronizable pairs of maps and show that synchronization may fail even when the images of the maps intersect. This leads to genuinely bivariate rigidity phenomena, including the appearance of infinite values and the lack of homotopy invariance, which have no analogue in the classical or univariate settings.

We then develop the structural theory of $\mathrm{TC}(f,g)$. This includes symmetry, behavior under pre- and post-composition, product inequalities, stability results, and several bounds relating $\mathrm{TC}(f,g)$ to the topological complexity of the target space. A central result in this direction is the collaboration principle, which shows that when one of the maps is a fibration, the complexity of synchronization is controlled by the other map. This principle explains how one subsystem may absorb the coordination cost and situates the bivariate invariant within the broader framework of sectional category and topological complexity of maps.

To complement the strict theory, we introduce a homotopy–invariant bivariate topological complexity $\mathrm{TC}_H(f,g)$ of Scott type. This invariant is defined in terms of homotopic distance and provides a robust lower bound for $\mathrm{TC}(f,g)$. We study its basic properties, its relationship with homotopic distance, and several computable bounds, highlighting both its flexibility and its limitations when compared with the strict invariant.

A systematic comparison between $\mathrm{TC}_H(f,g)$ and $\mathrm{TC}(f,g)$ reveals a hierarchy of rigidity. In favorable situations, such as when the maps are fibrations, both invariants coincide. In general, however, they may differ substantially, exhibiting strict gaps that reflect the cost of enforcing exact synchronization as opposed to homotopical compatibility.

Finally, we develop cohomological estimates for $\mathrm{TC}(f,g)$ and $\mathrm{TC}_H(f,g)$, extending classical zero–divisors techniques to the bivariate setting and providing computable obstructions to synchronization.

The paper is organized as follows. Section~1 recalls the necessary background on sectional invariants. Section~2 introduces Pave\v{s}i\'c-type bivariate topological complexity and its first properties. Section~3 develops cohomological estimates for $\mathrm{TC}(f,g)$. Section~4 is devoted to the case where one of the maps is a fibration, culminating in the collaboration principle. Section~5 introduces the Scott-type homotopy–invariant $\mathrm{TC}_H(f,g)$. Sections~6 and~7 compare the strict and homotopy–invariant theories and relate $\mathrm{TC}_H(f,g)$ to homotopic distance. Section~8 establishes further bounds for $\mathrm{TC}_H(f,g)$. We conclude with a discussion and outlook on future directions.

\section{Preliminaries}

Throughout this paper we shall assume that all the considered maps are continuous. In this section we recall the basic notions of sectional invariants that will underlie our constructions. Let $f:E \to B$ be a continuous map. The \emph{sectional category} or \emph{Svarc genus} of $f$, denoted $\mathrm{secat}(f)$, is the least integer $k$ (or infinity) such that $B$ admits an open cover $B = U_0 \cup \cdots \cup U_k$ with continuous maps $s_i:U_i \to E$ satisfying $f \circ s_i \simeq \mathrm{inc}_{U_i}$:
$$
\xymatrix{
{U_i} \ar@{^(->}[rr]^{\mathrm{inc}_{U_i}} \ar[dr]_{s_i} & & {B} \\
 & {E} \ar[ur]_f &  }
$$
in other words, each $U_i$ admits a local homotopy section of $f$. The (standard) \emph{sectional number} of $f$, denoted here by $\sec(f)$, is defined in the same way but requiring strict local sections, that is, $f \circ s_i = \mathrm{inc}_{U_i}$. Clearly one always has $\mathrm{secat}(f) \leq \mathrm{sec}(f)$, and both invariants coincide whenever $f$ is a fibration \cite{Schw}. 

A third variant, used by Pave\v{s}i\'c in \cite{Pav2} for the definition of the topological complexity of a map, is the filtered sectional number $\mathrm{sec}^*(f)$. This is the least integer $k$ for which there exists an increasing filtration of open subsets
\[
\emptyset = U_{-1} \subseteq U_0 \subseteq U_1 \subseteq \cdots \subseteq U_k = B,
\]
such that each difference $U_i \setminus U_{i-1}$ admits a strict local section of $f$. One always has $\mathrm{sec}^*(f) \leq \mathrm{sec}(f)$. The motivation for this definition is that, for maps with singularities, the version with open covers can become infinite, whereas the filtered definition often yields a finite invariant. On the other hand, when $f$ is a fibration between absolute neighborhood retracts (ANRs), all three invariants coincide  \cite{Pav2}. 

In this paper we shall consistently work with the sectional number $\mathrm{sec}(f)$, rather than the filtered number $\mathrm{sec}^*(f)$ adopted by Pave\v{s}i\'c in his univariate definition of topological complexity. This choice is deliberate, for several reasons. First, $\mathrm{sec}(f)$ is the formulation closest to Schwarz’s classical invariant, to Farber’s definition of topological complexity, and to the Lusternik–Schnirelmann category, all of which are traditionally expressed in terms of open covers rather than filtrations. Second, in the bivariant setting the unfiltered definition typically leads to more transparent proofs, while filtrations introduce an additional layer of rigidity without providing clear advantages for our purposes. Third, our interest here is not in dealing with singular behaviours of maps—where filtered versions are useful—but in highlighting structural phenomena that arise specifically in the bivariant framework; these are better captured by the strict, non-filtered formulation. We emphasize that in the case of fibrations—which constitute many of our main examples—the distinction disappears, as one always has $\mathrm{secat}(f)=\mathrm{sec}(f)$. Moreover, if the spaces involved are ANRs, then $\mathrm{sec}^*(f)=\mathrm{sec}(f)$ as well, so that all three invariants coincide in this setting.

\begin{remark} One could equally well develop the bivariant theory 
using the filtered invariant $\mathrm{sec}^*(f)$ instead of $\sec(f)$. 
Most of the basic results would then admit parallel statements, 
yielding a formally analogous theory. 
We leave this alternative approach to the interested reader.
\end{remark}

For later use we recall some classical notions that will appear throughout this work:

\begin{itemize}
\item The \emph{Lusternik–Schnirelmann category} of a space $X$, denoted $\mathrm{cat}(X)$, 
is the least integer $k$ such that $X$ can be covered by $k+1$ open subsets, each of which is contractible in $X$. When  $X$ is path-connected, this invariant admits an equivalent description in terms of sectional category: it coincides with the sectional category of the pointed path fibration
\[
p:PX \longrightarrow X,\qquad p(\gamma)=\gamma(1)
\]
where $PX$ denotes the space of paths in $X$ starting at the basepoint.

\item Given a continuous map $f:X\to Y$, the \emph{category of $f$}, denoted $\mathrm{cat}(f)$, 
is the least integer $k$ such that $X$ can be covered by $k+1$ open subsets 
on which the restriction of $f$ is null-homotopic.

\item The \emph{topological complexity} of a space $X$ in the sense of Farber \cite{Far}, 
denoted $\mathrm{TC}(X)$, is the sectional category of the free path fibration
\[
\pi_X:X^I \longrightarrow X\times X,\qquad \pi_X(\gamma)=(\gamma(0),\gamma(1)),
\]
where $X^I$ is the path space of $X$.

\item Following Pave\v{s}i\'c \cite{Pav2}, the \emph{topological complexity of a map} 
$f:X\to Y$, denoted $\mathrm{TC}(f)$, is defined as the filtered sectional number of the map
\[
\pi_f:X^I\longrightarrow X\times Y,
\]
where $\pi_f(\alpha )=(\alpha(0),f(\alpha (1)))$.
\end{itemize}

\begin{remark} 
In Pave\v{s}i\'c’s original definition the invariant $\mathrm{TC}(f)$ is based on the filtered 
sectional number $\mathrm{sec}^*$. In this paper, however, we shall adopt the unfiltered 
version $\mathrm{sec}$ when writing $\mathrm{TC}(f)$, both for simplicity and 
because it better suits the structural aspects of the bivariant setting. 
\end{remark}

\medskip
\noindent\textbf{Convention.} 
Throughout the paper we shall assume that the reader is familiar with these invariants 
and with their elementary properties, which we do not reproduce here in detail. 
Standard references are Farber’s foundational paper on topological complexity \cite{Far}, 
the monograph by Cornea–Lupton–Oprea–Tanré \cite{C-L-O-T}, and Pave\v{s}i\'c’s work 
on the topological complexity of maps \cite{Pav2}. 
Our focus will be on the new phenomena that arise in the bivariant setting.

\medskip
Having recalled these notions and fixed our conventions, we now turn to a few basic properties that will be used repeatedly in what follows. They describe the behavior of 
$\mathrm{secat}$ and $\mathrm{sec}$
 under (homotopy) commutative diagrams, pullbacks, and equivalences. The next three statements provide a convenient toolkit for the constructions developed later. The proofs of the first two results are straightforward and well-known. The last statement   can be found and proved in  \cite{IZ-T} (see also \cite{IZ} and \cite{IZ-G}.) We begin with the simplest case: comparing sectional invariants along a commutative or homotopy commutative triangle.

\begin{proposition}\label{triangle}
Consider the following diagram
$$
\xymatrix{
{X'} \ar[rr] \ar[dr]_{f'} & & {X} \ar[dl]^f \\
& Y & }$$
\begin{enumerate}
\item If the diagram is homotopy commutative, then $\mathrm{secat}(f)\leq \mathrm{secat}(f').$

\item If the diagram is strictly commutative, then $\mathrm{sec}(f)\leq \mathrm{sec}(f')$ (and $\mathrm{secat}(f)\leq \mathrm{secat}(f')$).
\end{enumerate}
\end{proposition}

Next we record the behavior of sectional invariants with respect to pullbacks. This property will be essential when working with fibered constructions in later sections.

\begin{proposition}\label{pullback}
Let
$$\xymatrix{
{W} \ar[r] \ar[d]_{f'} & {X} \ar[d]^f \\
{Z} \ar[r] & {Y} }$$
\noindent be a pullback. Then $\mathrm{secat}(f')\leq \mathrm{secat}(f)$ and  $\mathrm{sec}(f')\leq \mathrm{sec}(f)$. 
\end{proposition}

We also include a slightly more elaborate statement involving a square diagram, which relates sectional invariants across both vertical and horizontal maps. This result will allow us to transfer estimates through homotopy equivalences and homeomorphisms.

\begin{proposition}\label{invariance}

Consider the following diagram
$$\xymatrix{
{W} \ar[r]^{\varphi } \ar[d]_{f'} & {X} \ar[d]^f \\
{Z} \ar[r]_{\psi } & {Y} }$$
\begin{enumerate}
\item If the diagram is homotopy commutative, then
$$(\mathrm{secat}(f')+1)\cdot (\mathrm{secat}(\psi )+1)\geq \mathrm{secat}(f)+1.$$
If, in addition, $\varphi $ and $\psi $ are homotopy equivalences, then $\mathrm{secat}(f)=\mathrm{secat}(f').$ 

\item If the diagram is strictly commutative, then
$$(\mathrm{sec}(f')+1)\cdot (\mathrm{sec}(\psi )+1)\geq \mathrm{sec}(f)+1.$$ If, in addition, $\varphi $ and $\psi $ are homeomorphisms, then $\mathrm{sec}(f)=\mathrm{sec}(f')$ (and $\mathrm{secat}(f)=\mathrm{secat}(f')$). 
\end{enumerate}
\end{proposition}

\section{Bivariate topological complexity: definition and first properties}

Topological complexity was originally introduced by Farber to capture the intrinsic discontinuities of motion planning algorithms on a single configuration space. 
Later, Pave\v{s}i\'c extended this idea to the setting of a single map $f:X\to Z$, defining the topological complexity of a map in terms of strict sectional invariants. 
Our aim in this section is to go one step further and consider the simultaneous behavior of two maps with a common codomain. 
This leads to the notion of \emph{Pave\v{s}i\'c-type bivariate topological complexity} or simply \emph{bivariate topological complexity} when no ambiguity arises, which quantifies the minimal amount of local data required to coordinate two systems whose dynamics are constrained by a shared target space. 
In this framework, the corresponding motion planning problem involves pairs of paths that must remain synchronized through the maps into the common base.

Consider a cospan of maps $X \xrightarrow{f} Z \xleftarrow{g} Y$. The \emph{bivariate motion planning problem} consists of constructing an algorithm that, given as input a pair of points $(x, y) \in X \times Y$, produces as output a pair of continuous paths $\alpha \in X^I$, $\beta \in Y^I$ satisfying the conditions:
$$
f \circ \alpha = g \circ \beta, \quad \alpha(0) = x, \quad \beta(1) = y.
$$
Let us present this construction in a different way. 

\begin{definition} Given the pullback of the maps $f$ and $g$,
$$\Delta_{f,g} := \{(x,y) \in X \times Y \mid f(x) = g(y)\}.
$$
A pair of paths $(\alpha, \beta) \in X^I \times Y^I$ are $(f,g)$-synchronized if $f \circ \alpha = g \circ \beta$. So, the corresponding path-space is denoted by $\Delta_{f,g}^I$:
$$
\Delta_{f,g}^I = \{(\alpha, \beta) \in X^I \times Y^I \mid f \circ \alpha = g \circ \beta \},
$$
and introduce the $(f,g)$-endpoint evaluation map
$$
\pi_{f,g} : \Delta _{f,g}^I \to X \times Y, \quad \pi_{f,g}(\alpha, \beta) = (\alpha(0), \beta(1)).
$$

A \textit{bivariate motion planning algorithm} for the pair $(f, g)$ is a (not necessarily continuous) section $s : X \times Y \to \Delta _{f,g}^I$ of the map $\pi_{f,g}$, that is, $\pi_{f,g} \circ s = \mathrm{id}_{X \times Y}.$
This situation can be represented by the commutative diagram:
$$
\xymatrix{
X \times Y \ar[rr]^{\mathrm{id}_{X \times Y}} \ar[dr]_{s} & & X \times Y \\
& \Delta _{f,g}^I \ar[ur]_{\pi_{f,g}} &
}
$$

\end{definition}

The next result establishes how the existence of a global continuous section of $\pi_{f,g}$ and the homotopical triviallity of $f$ and $g$ are related.

\begin{proposition}\label{global}
If there exists a global continuous bivariate motion planning algorithm $s : X \times Y \to \Delta _{f,g}^I$, then the maps $f$ and $g$ must both be null-homotopic.    
\end{proposition}

\begin{proof}
If $s=(s_1,s_2)$, then $s_1(x,y)\in X^I,$ $s_2(x,y)\in Y^I$ are such that $f\circ s_1(x,y)=g\circ s_2(x,y)$, and $s_1(x,y)(0)=x$, $s_2(x,y)(1)=y,$ for all $(x,y)\in X\times Y.$  Now, fix a point $y_0\in Y$ and consider the homotopy
$H:X\times I\to Z$ given as $$H(x,t):=f(s_1(x,y_0)(t)).$$
Then, $H(x,0)=f(x)$ and $H(x,1)=g(y_0),$ proving that $f$ is null-homotopic. Similarly, $g$ is also null-homotopic.
\end{proof}

Since in most applications the maps $f$ and $g$ are not null-homotopic, we are led to consider the key notion of bivariate topological complexity of a pair of maps.

\begin{definition}
Let $X \stackrel{f}{\longrightarrow} Z \stackrel{g}{\longleftarrow} Y$ be continuous maps. The \emph{bivariate topological complexity} of the pair $(f, g)$ is defined as the sectional number of the map $\pi_{f,g}$:
$$
\mathrm{TC}(f,g) := \mathrm{sec}(\pi_{f,g}).
$$
In other words, $\mathrm{TC}(f,g)$ is the smallest integer $n \geq 0$ such that $X \times Y$ can be covered by $n+1$ open sets $\{U_i\}_{i=0}^n$, on each of which there exists a strict local section $s_i : U_i \to \Delta _{f,g}^I$ of $\pi_{f,g}$:
$$
\xymatrix{
U_i \ar@{^(->}[rr] \ar[dr]_{s_i} & & X \times Y \\
& \Delta_{f,g}^I \ar[ur]_{\pi_{f,g}} &
}
$$
If no such finite cover exists, we set $\mathrm{TC}(f,g) = \infty$.
\end{definition}

It is important to point out that the condition $f\circ \alpha = g\circ \beta$ means that the two paths 
$\alpha$ in $X$ and $\beta$ in $Y$ evolve in perfect synchrony 
through the target space $Z$: at every parameter $t\in I$ their images 
coincide, $f(\alpha(t))=g(\beta(t))$. Formally, $t$ is just the path parameter, 
but it can naturally be interpreted as a ``time'' variable. In this sense, 
$\mathrm{TC}(f,g)$ quantifies the complexity of enforcing a 
\emph{strict synchronization of trajectories in the target space}.

\medskip
Clearly, $\mathrm{TC}(id_X,id_X)=\mathrm{TC}(X)$ is the classical topological complexity. Moreover, note that $\pi_{f,g}$ is not a Hurewicz fibration in general. Consequently, its sectional number may differ from its sectional category. In particular, $\mathrm{TC}(f,g)$ is not, in general, a homotopy invariant.

\begin{remark}\label{ex-infty}
Proposition~\ref{global} shows that $\mathrm{TC}(f,g)=0$ implies both $f$ and $g$ are
null-homotopic. The converse, however, fails in general.
For instance, suppose $f=c_{z}$ and $g=c_{z'}$ are constant maps at two distinct points 
$z,z'\in Z$. In this case the pullback $\Delta_{f,g}$ is empty, so the map 
$\pi_{f,g}$ admits no local sections and hence $\mathrm{TC}(f,g)=\infty$.

This example also illustrates that $\mathrm{TC}(f,g)$ is not a homotopy invariant.
Indeed, if $Z$ is path-connected, then $g=c_{z'}$ is homotopic to $g'=c_{z}$. 
In this situation $\Delta_{f,g'}=X\times Y$ and there exists a global section of 
$\pi_{f,g'}$ given by constant paths $\alpha(t)\equiv x$, $\beta(t)\equiv y$. 
Therefore $\mathrm{TC}(f,g')=0$, while $\mathrm{TC}(f,g)=\infty$.
Thus, replacing $g$ by a homotopic map can drastically alter the value of the invariant.
\end{remark}

Now we check that bivariate topological complexity is, indeed, a generalization of Pavesic's \cite{Pav2} topological complexity of a map. 

\begin{proposition}
Let $f:X\to Z$ be a continuous map. Then
$$\mathrm{TC}(f)=\mathrm{TC}(f,id_Z).$$
\end{proposition}

\begin{proof}
There exists  $\Phi :\Delta _{f,id_Z}^I=\{(\alpha ,\beta )\in X^I\times Z^I \mid f\circ \alpha = \beta \}\to X^I$ a clear homeomorphism, given by $\Phi (\alpha ,\beta ):=\alpha .$ Moreover, this homeomorphism fits into the following commutative triangle:
$$
\xymatrix{
{\Delta _{f,id_Z}^I} \ar[rr]^{\Phi }_{\cong } \ar[dr]_{\pi _{f,id_Z}} & & {X^I} \ar[dl]^{\pi _f} \\
 & {X\times Z} & }$$
Therefore, by Proposition \ref{invariance} (2), we have that $\mathrm{TC}(f)=\mathrm{sec}(\pi _f)=\mathrm{sec}(\mathrm{\pi _{f,id_Z}})=\mathrm{TC}(f,id_Z).$
\end{proof}

\medskip
In practical scenarios it is essential to restrict attention to pairs of maps that can actually be synchronized. Let us formalize this notion.

\begin{definition}  Let $X \stackrel{f}{\longrightarrow} Z \stackrel{g}{\longleftarrow} Y$ be a cospan of maps. The pair $(f,g)$ is \emph{synchronizable} if, for every $x\in X$ and $y\in Y$, there exist paths $\alpha:I\to X$ and $\beta:I\to Y$ such that the pair $(\alpha ,\beta )$ is $(f,g)$-synchronized and
$$
\alpha(0)=x, \qquad \beta(1)=y.
$$  
In other words, any initial configuration in $X$ and final configuration in $Y$ can be connected through coordinated motions whose images under $f$ and $g$ remain perfectly matched in $Z$.

\end{definition}

The following characterization is straightforward:

\begin{lemma}
 Let $X \stackrel{f}{\longrightarrow} Z \stackrel{g}{\longleftarrow} Y$ be a cospan of maps. Then $(f,g)$ is synchronizable if and only if the map $\pi _{f,g}$ is surjective.
\end{lemma}

Consequently, if a pair $(f,g)$ fails to be synchronizable, then the associated bivariate topological complexity $\mathrm{TC}(f,g)$ is necessarily infinite. 

\begin{example}\label{ex-sync}

Let us present some illustrative examples of synchronizable and non-synchronizable pairs:

\begin{enumerate}

 \item Let $f:X\to Z$ be a continuous surjection with $X$ path-connected. Then the pair $(f,\mathrm{id}_Z)$ is always synchronizable. Indeed, given $x\in X$ and $z\in Z$, choose a path $\alpha:I\to X$ with $\alpha(0)=x$ and $\alpha(1)\in f^{-1}(z)\neq \emptyset$. Setting $\beta:=f\circ \alpha$, we have $f\circ \alpha = \mathrm{id}_Z\circ \beta$, with $\alpha(0)=x$ and $\beta(1)=z$.

\item Consider maps $X \stackrel{f}{\longrightarrow} Z \stackrel{g}{\longleftarrow} Y$, where $f$ is surjective, $X$ is path-connected, and $g$ is a surjective fibration. Then $(f,g)$ is synchronizable. In fact, given $x\in X$ and $y\in Y$, choose $x'\in X$ with $f(x')=g(y)$. By path-connectedness, there exists a path $\alpha:I\to X$ with $\alpha(0)=x$ and $\alpha(1)=x'$. Since $g$ is a fibration, we can lift $f\circ \alpha$ to a path $\beta:I\to Y$ with $\beta(1)=y$ in the diagram
$$
\xymatrix{
 {} & (Y,y) \ar[d]^g \\
 {(I,1)} \ar[r]_{f\circ \alpha} \ar@{.>}[ur]^{\beta} & {(Z,g(y))}.
}
$$
Thus $f\circ \alpha = g\circ \beta$, with $\alpha(0)=x$ and $\beta(1)=y$.

\item A trivial case where $(f,g)$ is not synchronizable occurs when $f(X)\cap g(Y)=\emptyset$. In this situation the pullback $\Delta_{f,g}$ is empty, so the synchronicity condition cannot be satisfied. Consequently, $\mathrm{TC}(f,g)=\infty$.

\item Let $Z=S^1\vee S^1$ be the wedge of two circles with common basepoint $*$. 
Take $X=Y=S^1$, and define $f,g\colon S^1\to Z$ as the inclusions into the first and second circle, respectively. 
Then $f(X)\cap g(Y)=\{*\}\neq\emptyset$. 
Nevertheless, $(f,g)$ is not synchronizable. Pick $x\in X$ with $f(x)\neq *$ and $y\in Y$ with $g(y)\neq *$. 
If there existed paths $\alpha:I\to X$ and $\beta:I\to Y$ with $\alpha(0)=x$, $\beta(1)=y$ and 
$f\circ \alpha \;=\; g\circ \beta,$
their common image would be contained in $f(X)\cap g(Y)=\{*\}$, forcing $f(x)=*$, a contradiction. 
Hence $(f,g)$ is not synchronizable, and consequently
$\mathrm{TC}(f,g)=\infty.$
\end{enumerate}
\end{example}

More generally, synchronizability can be inherited from factorizations through a common subspace:

\begin{proposition}\label{prop:reduction-common-W}
Let \(i:W\hookrightarrow Z\) be the inclusion of a path–connected subspace. Suppose
\[
f=i\circ f'\colon X\longrightarrow Z,\qquad g=i\circ g'\colon Y\longrightarrow Z
\]
for some maps \(f'\colon X\to W\) and \(g'\colon Y\to W\).
If the pair \((f',g')\) is synchronizable, then the pair \((f,g)\) is synchronizable.
\end{proposition}

\begin{proof}
By hypothesis, for every \(x\in X\) and \(y\in Y\) there exist paths \(\alpha\colon I\to X\) and \(\beta\colon I\to Y\) such that
$f'\circ \alpha ,$ $ g'\circ \beta,$ $\alpha(0)=x,\beta(1)=y.$
Composing with the inclusion \(i\colon W\hookrightarrow Z\) yields
\[
f\circ \alpha=(i\circ f')\circ \alpha=i\circ (f'\circ \alpha)
=i\circ (g'\circ \beta)=(i\circ g')\circ \beta=g\circ \beta,
\]
with \(\alpha(0)=x\) and \(\beta(1)=y\). Hence \((f,g)\) is synchronizable in \(Z\).
\end{proof}

\begin{corollary}\label{cor:fibration-pattern}
With the same notation as in Proposition~\ref{prop:reduction-common-W}, assume that \(W\) is path–connected. Then \((f,g)\) is synchronizable in each of the following cases:
\begin{enumerate}[label=(\alph*)]
\item \(f'\) is surjective and \(g'\) is a fibration (for instance, a covering projection).
\item \(g'\) is surjective and \(f'\) is a fibration.
\end{enumerate}
\end{corollary}

\begin{proof}
Just use Proposition~\ref{prop:reduction-common-W} above and Example \ref{ex-sync} (2).
\end{proof}

\begin{remark}
This criterion yields synchronizable pairs where neither \(f\) nor \(g\) is surjective or a fibration onto \(Z\): it suffices to factor them through a common subspace \(W\) where synchronizability can be established.  
\item For instance, let \(Z=\mathbb C\), \(W=S^1\subset \mathbb C\), and define \(f,g\colon S^1\to Z\) by \(f(z)=z^2\), \(g(z)=-z^2\). Writing \(f=i\circ q\) and \(g=i\circ (\mu_{-1}\circ q)\) with \(q\colon S^1\to S^1\) the degree–two covering and \(\mu_{-1}(z)=-z\), the pair \((q,\mu_{-1}\circ q)\) is synchronizable in \(W\) by the path–lifting property of coverings. Proposition~\ref{prop:reduction-common-W} then implies that \((f,g)\) is synchronizable in \(Z\).
\end{remark}

Now we present a concrete example that illustrates a computation of the bivariate topological complexity.

\begin{example}\label{iconic-example}
Let $f,g:S^1 \to \mathbb{C}$ be defined by $f(z)=z^2$ and $g(z)=-z^2$ for all $z\in S^1$. Then we have $\mathrm{TC}(f,g)=1$. 
Indeed, a pair of paths $(\alpha,\beta)\in \Delta^I_{f,g}$ must satisfy $\alpha(t)^2=-\beta(t)^2$, and therefore, $\alpha(t)/\beta(t)\in \{-i,i\}$ for all $t\in I$. By continuity and the connectedness of $I$, the quotient $\alpha(t)/\beta(t)$ is constant, hence either $\alpha(t)=i\beta(t)$ for all $t$, or $\alpha(t)=-i\beta(t)$ for all $t$. Consequently, $\Delta^I_{f,g}$ splits into two path components,
\[
\Delta^I_{f,g}=E_+ \sqcup E_-,
\]
with 
$E_+=\{(\alpha,\beta)\in (S^1)^I\times (S^1)^I \mid \alpha=i\beta\}$
and 
$E_-=\{(\alpha,\beta)\in (S^1)^I\times (S^1)^I \mid \alpha=-i\beta\}.$

For each component we obtain a commutative diagram:
\[
\xymatrix{
(S^1)^I \ar[d]_{\pi_{S^1}} \ar[rr]^{\Phi_\pm}_{\cong} & & E_\pm \ar[d]^{\pi_{f,g}|_{E_\pm}} \\
S^1\times S^1 \ar[rr]^{\psi_\pm}_{\cong} & & S^1\times S^1
}
\]
where $\pi_{S^1}$ denotes the standard path fibration, $\Phi_+(\beta)=(i\beta,\beta)$, $\Phi_-(\beta)=(-i\beta,\beta)$, and $\psi_\pm(u,v)=(\pm i u, v)$. The horizontal maps are homeomorphisms, hence 
\[
\mathrm{sec}(\pi_{f,g}|_{E_\pm})=\mathrm{sec}(\pi_{S^1})=\mathrm{TC}(S^1)=1.
\]

First, let us show that $\mathrm{TC}(f,g)\geq 1$.  Suppose, for contradiction, that $\mathrm{TC}(f,g)=0$. Then $\pi_{f,g}$ would admit a global strict section $\sigma:S^1\times S^1\to \Delta^I_{f,g}$. Since $S^1\times S^1$ is path-connected, the image of $\sigma$ would lie entirely in one component, say $E_+$ or $E_-$. This would imply $\mathrm{sec}(\pi_{f,g}|_{E_\pm})=0$, i.e.~$\mathrm{TC}(S^1)=0$, a contradiction. Therefore $\mathrm{TC}(f,g)\geq 1$. 

To show that $\mathrm{TC}(f,g)\leq 1$, we construct local strict sections over a two-set open cover of $S^1\times S^1$. Consider the open subsets:
\[
U_+:=\{(u,v)\in S^1\times S^1 \mid -\frac{v}{iu}\neq 1\}; \hspace{8pt} U_{-}:=\{(u,v)\in S^1\times S^1 \mid \frac{v}{iu}\neq 1\}
\]
Since every point $(u,v)$ lies in at least one of these sets, we have $U_+\cup U_{-}=S^1\times S^1.$
Choose a continuous branch of the argument $\mathrm{Arg}_{1}:S^1\setminus\{1\}\to (0,2\pi),$ and define 
$$\theta (u,v):=\mathrm{Arg}_1(-\frac{v}{iu}),\hspace{6pt}\mbox{for all}\hspace{4pt}(u,v)\in U_+,$$
\noindent and 
$$\begin{cases} s_2(u,v)(t):=-i\,u\, e^{\,i\,\theta (u,v)\,t}, 
\\
s_1(u,v)(t):=i\,s_2(u,v)(t)=u\,e^{\,i\,\theta (u,v)\,t}. \end{cases}
$$
Then $s=(s_1,s_2):U_+\to E_+\subseteq \Delta^I_{f,g}$ is a strict local section of $\pi_{f,g}$. 
Similarly, define $\theta '(u,v):=\mathrm{Arg}_1(\frac{v}{iu}),$ for all $(u,v)\in U_{-}$ and
$$
\begin{cases} s'_2(u,v)(t):=i\,u e^{\,i\,\theta '(u,v)\,t}, \\
s'_1(u,v)(t):=-i\,s'_2(u,v)(t)=u\,e^{\,i\,\theta '(u,v)\,t}. \end{cases}
$$
Hence $s'=(s'_1,s'_2):U_-\to E_-\subseteq \Delta^I_{f,g}$ is another strict local section of $\pi _{f,g}$. We conclude that $\mathrm{TC}(f,g)=1$.
\end{example}

Our first structural property establishes that the bivariate topological complexity is symmetric with respect to its two arguments. 

\begin{proposition} Let $X \stackrel{f}{\longrightarrow} Z \stackrel{g}{\longleftarrow} Y$ be a cospan of maps. Then we have that $$\mathrm{TC}(f,g)=\mathrm{TC}(g,f).$$
\end{proposition}

\begin{proof}
We have just to take into account the following commutative diagram
$$\xymatrix{
{\Delta _{f,g}^I} \ar[d]_{\pi _{f,g}} \ar[rr]^{T}_{\cong } & & {\Delta _{g,f}^I} \ar[d]^{\pi _{g,f}} \\
{X\times Y} \ar[rr]_t^{\cong } & & {Y\times X} }$$
\noindent where the homeomorphisms $T$ and $t$ are respectively defined as $T(\alpha ,\beta ):=(\bar{\beta },\bar{\alpha })$ and $t(x,y):=(y,x).$ Here $\bar{\alpha }$ and $\bar{\beta }$ stand for the inverse paths of $\alpha $ and $\beta .$ Consequently, using Proposition 1.5 we have
$$\mathrm{TC}(f,g)=\mathrm{sec}(\pi _{f,g})=\mathrm{sec}(\pi _{g,f})=\mathrm{TC}(g,f)$$
\end{proof}

Our second result establishes an interesting inequality for general maps with a common target space. 
\begin{proposition}\label{bound}
Let \(X \xrightarrow{\,f\,} Z \xleftarrow{\,g\,} Y\) be a cospan of maps. The following inequality holds.

$$(\mathrm{TC}(f,g)+1)\cdot (\mathrm{sec}(f\times g)+1)\geq \mathrm{TC}(Z)+1.$$
\end{proposition}

\begin{proof}

We consider the following strictly commutative square
$$\xymatrix{
{\Delta ^I_{f,g}} \ar[d]_{\pi _{f,g}} \ar[rr]^{\rho } & & {Z^I} \ar[d]^{\pi _Z} \\
{X\times Y} \ar[rr]_{f\times g} & & {Z\times Z}
}$$
\noindent where $\rho (\alpha ,\beta ):=f\circ \alpha \hspace{3pt}(=g\circ \beta ).$, for all $(\alpha ,\beta )\in \Delta ^I_{f,g}.$ By means of \cite[Prop. 2.3 (b)]{IZ-T}, we obtain that
$$ (\mathrm{sec}(\pi _{f,g})+1)\cdot (\mathrm{sec}(f\times g)+1)\geq \mathrm{sec}(\pi _Z)+1.$$
\noindent and thus, by means of the definitions of $\mathrm{TC}(f,g)$ and $\mathrm{TC}(Z)$, the desired inequality follows. 

\end{proof}

\begin{corollary}
Let \(X \xrightarrow{\,f\,} Z \xleftarrow{\,g\,} Y\) be a cospan of maps where $f$ and $g$ admit strict global sections. Then

$$ \mathrm{TC}(f,g)\geq \mathrm{TC}(Z). $$
\end{corollary}

\begin{proof}

Since $f$ and $g$ admit global sections, then $f\times g$ also admits a global section and thus, $\mathrm{sec}(f\times g)=0$. Now, taking into account the previous proposition, we conclude the result.
  
\end{proof}

Another interesting structural feature of bivariate topological complexity emerges when considering products of maps. 
The following result was proved in \cite[Prop. 22, p.84]{Schw}. We agree that a normal space is, by definition, required to be Hausdorff.

\begin{proposition}\label{lem:sec-product}
Let $f:X\to Y$ and $f':X'\to Y'$ be continuous maps, and assume that $Y\times Y'$ is normal. Then
\[
\sec(f\times f')\ \leq\ \sec(f)+\sec(f').
\]
\end{proposition}

\begin{corollary}
Let $Z\times Z$ be a normal space (e.g. when $Z$ is metrizable) and let \(X \xrightarrow{\,f\,} Z \xleftarrow{\,g\,} Y\) be a cospan of maps. Then
$$ (\mathrm{TC}(f,g)+1)\cdot (\mathrm{sec}(f) + \mathrm{sec}(g) +1) \geq \mathrm{TC}(Z)+1. $$
If, in addition, $g$ admits a global section, then
$$ (\mathrm{TC}(f,g)+1)\cdot (\mathrm{sec}(f) +1) \geq \mathrm{TC}(Z) +1,$$
and similarly when $f$ admits a global section.
\end{corollary}

\begin{proof}
Just consider Proposition \ref{bound} and Proposition \ref{lem:sec-product}.

\end{proof}

Using Proposition \ref{lem:sec-product} above, we also obtain the subadditivity of bivariate topological complexity.

\begin{theorem}\label{thm:subadditivity-Pav}
Let
$
X \xrightarrow{\,f\,} Z \xleftarrow{\,g\,} Y, 
$ and $X' \xrightarrow{\,f'\,} Z' \xleftarrow{\,g'\,} Y'
$
be cospans of continuous maps, and suppose that the product $X\times Y\times X'\times Y'$ is normal 
(for instance, when these spaces are metrizable). Then
\[
\TC(f\times f',\,g\times g')\ \leq\ \TC(f,g)+\TC(f',g').
\]
\end{theorem}

\begin{proof}
Consider the commutative diagram
\[
\xymatrix{
{\Delta ^I_{f\times f',g\times g'}} \ar[d]_{\pi _{f\times f',g\times g'}} \ar[rr]^{T}_{\cong } & & {\Delta ^I_{f,g}\times \Delta ^I_{f',g'}} \ar[d]^{\pi _{f,g}\times \pi _{f',g'}} \\
{X\times X'\times Y\times Y'} \ar[rr]_{t}^{\cong } & & {X\times Y\times X'\times Y'}
}
\]
with the homeomorphisms
$
T\big((\alpha,\alpha'),(\beta,\beta')\big)=\big((\alpha,\beta),(\alpha',\beta')\big),$ and $
t(x,x',y,y')=(x,y,x',y').$
Hence
\[
\begin{aligned}
\TC(f\times f',\,g\times g') 
&= \sec(\pi _{f\times f',g\times g'}) \\
&= \sec(\pi _{f,g}\times \pi _{f',g'}) \\
&\leq \sec(\pi _{f,g})+\sec(\pi _{f',g'}) \\
&= \TC(f,g)+\TC(f',g'),
\end{aligned}
\]
by Proposition \ref{invariance} (2) and Proposition~\ref{lem:sec-product}.
\end{proof}

The next result establishes that TC's of pairs of maps are lower bounds for the TC of their corresponding pair of product maps.

\begin{proposition}[Lower bound by factors]\label{prop:lower-max}
Let
$
X \xrightarrow{\,f\,} Z \xleftarrow{\,g\,} Y$
and $
X' \xrightarrow{\,f'\,} Z' \xleftarrow{\,g'\,} Y'
$
be cospans of continuous maps. Then
\[
\max\{\TC(f,g),\,\TC(f',g')\}\ \le\ \TC(f\times f',\,g\times g').
\]
\end{proposition}

\begin{proof}
Consider the strictly commutative diagram
$$\xymatrix{
{\Delta ^I_{f\times f',g\times g'}} \ar[d]_{\pi _{f\times f',g\times g'}} \ar[rr]^{\varphi } & & {\Delta ^I_{f,g}} \ar[d]^{\pi _{f,g}} \\
{X\times X'\times Y\times Y'} \ar[rr]_(.6){\psi } & & {X\times Y} }$$
\noindent where $\varphi ((\alpha ,\alpha '),(\beta ,\beta ')):=(\alpha ,\beta )$ and $\psi (x,x',y,y')=(x,y),$ respectively. Choose points $x'_0\in X'$ and $y'_0\in Y'$ and define a global section of $\psi $ by
$$\iota :X\times Y\to X\times X'\times Y\times Y',\qquad \iota (x,y)=(x,x'_0,y,y'_0).$$
This proves that $\mathrm{sec}(\psi )=0.$ By Proposition \ref{invariance} (2) we have that $\mathrm{TC}(f,g)\leq \mathrm{TC}(f\times f',g\times g').$ 

%Fix basepoints $x'_0\in X'$ and $y'_0\in Y'$ and consider 
%$\iota:X\times Y\to X\times X'\times Y\times Y'$, $\iota(x,y)=(x,x'_0,y,y'_0)$. 
%As in the proof of Theorem~\ref{thm:subadditivity-Pav}, one checks that 
%$\iota^{*}(\Delta^I_{f\times f',\,g\times g'})\cong \Delta^I_{f,g}$ and the projection is $\pi_{f,g}$. {\red Since [Ref]} strict sectional category does not increase under pullback, hence 
%$\TC(f,g)\le \TC(f\times f',\,g\times g')$. 
Interchanging the roles of $(f,g)$ and $(f',g')$ gives the other inequality; taking the maximum yields the result. 
\end{proof}

\begin{corollary}
Under the hypotheses of Theorem~\ref{thm:subadditivity-Pav}, we have 
\[
\max\{\TC(f,g),\,\TC(f',g')\}\ \le\ \TC(f\times f',\,g\times g')\ \le\ \TC(f,g)+\TC(f',g').
\]
\end{corollary}

\begin{corollary}\label{cor:stability-zero}
Under the hypotheses of Theorem~\ref{thm:subadditivity-Pav}, if moreover $\TC(f',g')=0$, then
\[
\TC(f\times f',\,g\times g')=\TC(f,g).
\]
\end{corollary}

\begin{proof}
By Proposition~\ref{prop:lower-max}, $\TC(f,g)\le \TC(f\times f',\,g\times g')$. 
By Theorem~\ref{thm:subadditivity-Pav}, $\TC(f\times f',\,g\times g')\le \TC(f,g)+\TC(f',g')=\TC(f,g)$. 
Thus equality holds.
\end{proof}

\begin{corollary}
Let $X \xrightarrow{\,f\,} Z \xleftarrow{\,g\,} Y$ be a cospan of continuous maps, and let $K$ be any topological space. 
Assume that $X\times Y\times K\times K$ is normal (e.g.\ when $X$, $Y$, and $K$ are metrizable). Then
\[
\TC(f\times \mathrm{id}_{K},\, g\times \mathrm{id}_{K}) \;\leq\; \TC(f,g)\;+\;\TC(K).
\]
If, in addition, $K$ is contractible, then $\TC(f\times \mathrm{id}_{K},\, g\times \mathrm{id}_{K})=\TC(f,g).$
\end{corollary}

To conclude this section, we investigate how the complexity of a pair \((f,g)\) behaves under pre- and post-composition. We begin with the effect of post-composing with a map.  

\begin{proposition}
Let \(X \xrightarrow{\,f\,} Z \xleftarrow{\,g\,} Y\) be a cospan of maps, and let \(w:Z\to Z'\) be any map. Then
$$\mathrm{TC}(w\circ f,w\circ g)\leq \mathrm{TC}(f,g).
$$
Moreover, if there exists a map \(r:Z'\to Z\) such that \(r\circ w\circ f=f\) and \(r\circ w\circ g=g\) (for instance, if \(r\circ w=\mathrm{id}_Z\)), then
\[
\mathrm{TC}(w\circ f,w\circ g)=\mathrm{TC}(f,g).
\]
\end{proposition}

\begin{proof}
Every pair \((\alpha,\beta)\in \Delta^I_{f,g}\) also belongs to \(\Delta^I_{w\circ f,w\circ g}\), so there is a canonical map \(\Delta^I_{f,g}\to \Delta^I_{w\circ f,w\circ g}\) making the following diagram commute:
\[
\xymatrix{
{\Delta^I_{f,g}} \ar[rr] \ar[dr]_{\pi_{f,g}} & & {\Delta^I_{w\circ f,w\circ g}} \ar[dl]^{\pi_{w\circ f,w\circ g}} \\ & {X\times Y} &
}
\]
Thus, by Proposition 1.3 \(\mathrm{TC}(w\circ f,w\circ g)=\mathrm{sec}(\pi_{w\circ f,w\circ g})\leq \mathrm{sec}(\pi_{f,g})=\mathrm{TC}(f,g).\)

Now assume there exists \(r:Z'\to Z\) with \(r\circ w\circ f=f\) and \(r\circ w\circ g=g\). If \((\alpha,\beta)\in \Delta^I_{w\circ f,w\circ g}\), then \(w\circ f\circ \alpha = w\circ g\circ \beta\), and hence
\[
f\circ \alpha = r\circ w\circ f\circ \alpha = r\circ w\circ g\circ \beta = g\circ \beta.
\]
Therefore, there is a canonical map \(\Delta^I_{w\circ f,w\circ g}\to \Delta^I_{f,g}\) fitting into the commutative diagram
\[
\xymatrix{
{\Delta^I_{w\circ f,w\circ g}} \ar[rr] \ar[dr]_{\pi_{w\circ f,w\circ g}} & & {\Delta^I_{f,g}} \ar[dl]^{\pi_{f,g}} \\ & {X\times Y} &
}
\]
and consequently, applying again Proposition 1.3 we conclude that \(\mathrm{TC}(f,g)\leq \mathrm{TC}(w\circ f,w\circ g)\).
\end{proof}

We next analyze the effect of pre-composition. In this case, the behavior of the invariant depends on the properties of the maps used for the pre-composition.

\begin{proposition}\label{pre-comp-Pav}
Let \(X \xrightarrow{\,f\,} Z \xleftarrow{\,g\,} Y\) be a cospan of maps, and let \(u:X'\to X\), \(v:Y'\to Y\) be maps.
\begin{enumerate}

\item The following inequality holds:
$$(\mathrm{TC}(f\circ u,g\circ v)+1)\cdot (\mathrm{sec}(u\times v)+1)\geq \mathrm{TC}(f,g)+1.$$ Therefore, if \(u\) and \(v\) admit strict sections, then
$\mathrm{TC}(f,g)\leq \mathrm{TC}(f\circ u,g\circ v).$

\item If \(u\) and \(v\) are fibrations, then
$\mathrm{TC}(f\circ u,g\circ v)\leq \mathrm{TC}(f,g).
$

\item If \(u\) and \(v\) are fibrations admitting (homotopy) sections, then
$\mathrm{TC}(f\circ u,g\circ v)=\mathrm{TC}(f,g).
$
\end{enumerate}
\end{proposition}

\begin{proof}

(1) We consider the following strictly commutative square
$$\xymatrix{
{\Delta ^I_{f\circ u,g\circ v}} \ar[d]_{\pi _{f\circ u,g\circ v}} \ar[rr]^{\varphi } & & {\Delta ^I_{f,g}} \ar[d]^{\pi _{f,g}} \\
{X'\times Y'} \ar[rr]_{u\times v} & & {X\times Y}
}$$
\noindent where $\varphi (\alpha ,\beta ):=(u\circ \alpha ,v\circ \beta )$, for all $(\alpha ,\beta )\in \Delta ^I_{f\circ u,g\circ v}.$ Applying Proposition \ref{invariance} (2) we obtain the inequality. The rest of the proof is straightforward.

(2) Suppose \(u\) and \(v\) are fibrations. Let \(U\subseteq X\times Y\) be an open subset together with a strict local section
\[
\xymatrix{
{U} \ar@{^(->}[rr] \ar[dr]_{\sigma =(\sigma _1,\sigma _2)} & & {X\times Y} \\
& {\Delta ^I_{f,g}} \ar[ur]_{\pi _{f,g}} &
}
\]
Define \(V:=(u\times v)^{-1}(U)\subseteq X'\times Y'\). Consider the commutative diagram of solid arrows
\[
\xymatrix{
{V\times \{0\}} \ar[rr]^h \ar@{^(->}[d] & & {X'\times Y'} \ar[d]^{u\times v} \\
{V\times I} \ar[rr]_H \ar@{.>}[urr]^{\widetilde{H}} & & {X\times Y} }
\]
where \(h(x',y',0)=(x',y')\), and
\[
H(x',y',t)=(\sigma _1(u(x'),v(y'))(t),\;\sigma _2(u(x'),v(y'))(1-t)).
\]
Since \(u\times v\) is a fibration, the homotopy lifting property provides a lift \(\widetilde{H}=(\widetilde{H}_1,\widetilde{H}_2)\). Define
\[
\widetilde{\sigma }_1(x',y')(t):=\widetilde{H}_1(x',y',t), 
\quad
\widetilde{\sigma }_2(x',y')(t):=\widetilde{H}_2(x',y',1-t).
\]
Then \(\widetilde{\sigma}=(\widetilde{\sigma}_1,\widetilde{\sigma}_2):V\to \Delta^I_{f\circ u,g\circ v}\) is a strict local section of \(\pi_{f\circ u,g\circ v}\), proving the inequality.

(3) If \(u\) and \(v\) are fibrations admitting homotopy sections, the inequalities in (1) and (2) combine to yield the claimed equality.                    
\end{proof}

\section{Cohomological estimate for $\mathrm{TC}(f,g)$} 
 
To obtain general lower bounds valid for arbitrary pairs of continuous 
maps, we now turn to cohomological techniques. 
In this section we develop a bivariate analogue of the classical 
zero-divisors cup-length method, yielding a systematic estimate for 
$\mathrm{TC}(f,g)$ that requires no additional assumptions and extends 
the standard cohomological bound for Farber’s topological complexity.

Observe that, since in general $\mathrm{TC}(f,g)=\mathrm{sec}(\pi _{f,g})\geq \mathrm{secat}(\pi _{f,g})$, we have the natural cohomological lower bound:
$$\mathrm{nil}\hspace{3pt}\mathrm{Ker} (\pi _{f,g}^*)\leq \mathrm{TC}(f,g)$$
One can also take the inclusion $inc:\Delta _{f,g}\hookrightarrow X\times Y$ instead and consider $\mathrm{nil}\hspace{3pt}\mathrm{Ker} (inc^*).$ Nevertheless, we are interested in a more tractable estimate, even at the expense of some loss in sharpness. Given a cospan of maps \(X \xrightarrow{\,f\,} Z \xleftarrow{\,g\,} Y\), let us consider $I_Z$ the ideal of zero divisors
$$\mathrm{I}_Z:=\mathrm{Ker}(\Delta _Z^*)\subseteq H^*(Z\times Z)$$ \noindent where $\Delta _Z:Z\to Z\times Z$ denotes the diagonal map. Also consider the homomorphism $(f\times g)^*$ induced in cohomology by the product map $f\times g:X\times Y\to Z\times Z$. Then we define the image
$$\mathrm{I}_{f,g}:=(f\times g)^*(I_Z)$$
In general, this is not an ideal of $H^*(X\times Y)$. However, one can see that $\mathrm{I}_{f,g}$ is contained in $\mathrm{Ker} (\pi _{f,g}^*)$. In order to prove it show this, take the following commutative diagram of spaces and maps
$$\xymatrix{
{\Delta _{f,g}^I} \ar[rr]^{\rho } \ar[d]_{\pi _{f,g}} & & {Z^I} \ar[d]^{\pi _Z} & {Z} \ar[l]_c^{\simeq } \ar[dl]^{\Delta _Z} \\
{X\times Y} \ar[rr]_{f\times g} & & {Z\times Z} & 
}$$ \noindent where $\rho $ is defined as $\rho (\alpha ,\beta ):=f\circ \alpha \hspace{3pt}(=g\circ \beta ).$ This diagram induces another commutative diagram in cohomology:
$$\xymatrix{
 & {H^*(Z\times Z)} \ar[ld]_{\Delta _Z^*} \ar[d]_{\pi _Z^*} \ar[rr]^{(f\times g)^*} & & {H^*(X\times Y)} \ar[d]^{\pi _{f,g}^*} \\
 {H^*(Z)} & {H^*(Z^I)} \ar[l]^{c^*}_{\cong } \ar[rr]_{\rho ^*} & & {H^*(\Delta _{f,g}^I)} }$$
\noindent which proves that $\mathrm{I}_{f,g}=(f\times g)^*(\mathrm{I}_Z)\subseteq \mathrm{Ker} (\pi _{f,g}^*).$ 

We define $\mathrm{l.c.p.}\hspace{2pt}(\mathrm{I}_{f,g})$ as the least integer $k$ such that any product $u_0\smile \cdots \smile u_k$ of elements of $\mathrm{I}_{f,g}$ is null in $H^*(X\times Y)$ (here, $l.c.p.$ stands for the length of the cup product). Since $\mathrm{I}_{f,g}\subseteq \mathrm{Ker} (\pi _{f,g}^*)$, we have that $\mathrm{l.c.p.}\hspace{3pt}(\mathrm{I}_{f,g})\leq \mathrm{nil}\hspace{3pt}\mathrm{Ker} (\pi _{f,g}^*).$ 

\medskip
All this reasoning can be summarized in the following result:

\begin{theorem}\label{cohomological-bound}
For any cospan of continuous maps $X \stackrel{f}{\longrightarrow} Z \stackrel{g}{\longleftarrow} Y$ and any multiplicative cohomology theory $H^*$ we have the estimate
$$\mathrm{l.c.p.}\hspace{2pt}(\mathrm{I}_{f,g})\leq \mathrm{TC}(f,g).$$
\end{theorem}

\begin{remark}
In contrast with the classical case of the zero-divisor ideal, the subset $I_{f,g}$ is not necessarily an ideal of $H^*(X\times Y)$. 
Therefore the notion of nilpotency does not apply in this setting. The length of the cup product, however, is well defined for arbitrary subsets, and provides the appropriate cohomological estimate.
\end{remark}

For practical purposes we work with singular cohomology with coefficients in a field $\mathbb K$. 
By the Künneth theorem there is a natural isomorphism of graded $\mathbb K$-algebras 
$H^*(Z\times Z)\cong H^*(Z)\otimes H^*(Z)$, induced by the cross product 
$a\otimes b \mapsto pr_1^*(a)\smile pr_2^*(b)$. Under this identification, the map 
$\Delta_Z^*:H^*(Z\times Z)\to H^*(Z)$ agrees with the cup product 
$\smile:H^*(Z)\otimes H^*(Z)\to H^*(Z)$, so that 
$I_Z=\ker(\Delta_Z^*)=\langle\,u\otimes 1-1\otimes u\mid u\in H^*(Z)\,\rangle$.

For continuous maps $f:X\to Z$ and $g:Y\to Z$, the naturality of the cross product gives the commutative diagram
\[
\xymatrix{
 & & H^*(Z) \\
 H^*(Z)\otimes H^*(Z) \ar[urr]^{\smile}
    \ar[d]_{f^*\otimes g^*}
    \ar[rr]^{\times}_{\cong} 
  & & H^*(Z\times Z) \ar[u]_{\Delta_Z^*} \ar[d]^{(f\times g)^*} \\
 H^*(X)\otimes H^*(Y) \ar[rr]^{\times}_{\cong} & & H^*(X\times Y)
}
\]
and hence 
$I_{f,g}=(f\times g)^*(I_Z)=(f^*\otimes g^*)(I_Z)\subseteq H^*(X)\otimes H^*(Y)$. 
In particular, the images of the generators $u\otimes 1-1\otimes u$ are 
$f^*(u)\otimes 1-1\otimes g^*(u)$, which generate $I_{f,g}$.

\begin{remark}
Not every element of $I_Z$ is literally of the form $u\otimes 1-1\otimes u$. 
These elements generate the ideal, so more complicated terms (for instance $u\otimes u$ 
when $u\smile u=0$) can be written as products or linear combinations of them. 
\end{remark}

\begin{example}[Degrees $2$ and $3$ on the $2$-sphere]
Work over a field $\mathbb K$ with $\mathrm{char}\,\mathbb K=0$ (for instance, the field of rational numbers $\mathbb{Q}$). Let $X=Y=Z=S^2$ and pick continuous maps
$f,g:S^2\to S^2$ of degrees $2$ and $3$, respectively. 

Write $H^*(S^2)=\langle u\rangle$ with $|u|=2$ and $u^2=0$. Under the Künneth isomorphism
$H^*(Z\times Z)\cong H^*(Z)\otimes H^*(Z)$, the zero-divisor ideal is generated by $u\otimes 1-1\otimes u$, so
\[
I_{f,g}=(f^*\otimes g^*)(I_Z)
=\big\langle A\big\rangle,\qquad
A:=f^*(u)\otimes 1-1\otimes g^*(u)=2u\otimes 1-1\otimes 3u.
\]
Then $A\neq 0$ and, using the Koszul rule,
\[
A^2=(2u\otimes 1-1\otimes 3u)^2
=-(2u\otimes 1)(1\otimes 3u)-(1\otimes 3u)(2u\otimes 1)
=-12\,u\otimes u\neq 0,
\]
while $A^3=0$ for dimensional reasons. Hence
$
\mathrm{l.c.p.}\big(I_{f,g}\big)=2$ and, therefore,
$\mathrm{TC}(f,g)\geq 2.$

\end{example}

\begin{example}[Mixing different powers on the $5$-torus]
Work over a field $\mathbb K$ of characteristic $0$ (e.g.\ $\mathbb Q$). 
Let $Z=X=Y=T^5=(S^1)^5$ and write $H^*(T^5)=\Lambda(u_1,\dots,u_5)$ with $|u_i|=1$, 
$u_i\smile u_i=0$, and $u_1\smile\cdots\smile u_5$ generating $H^5(T^5)$.

Define continuous maps coordinatewise by
\[
f(z_1,\dots,z_5)=(z_1^{2},\,z_2^{3},\,z_3^{2},\,z_4^{4},\,z_5),\qquad
g(z_1,\dots,z_5)=(z_1,\,z_2^{2},\,z_3^{3},\,z_4,\,z_5^{4}).
\]
Then $f^*(u_1)=2u_1$, $f^*(u_2)=3u_2$, $f^*(u_3)=2u_3$, $f^*(u_4)=4u_4$, $f^*(u_5)=u_5$, and
$g^*(u_1)=u_1$, $g^*(u_2)=2u_2$, $g^*(u_3)=3u_3$, $g^*(u_4)=u_4$, $g^*(u_5)=4u_5$.
Under $H^*(Z\times Z)\cong H^*(Z)\otimes  H^*(Z)$ the zero-divisor ideal is generated by 
$u_i\otimes 1-1\otimes u_i$ ($i=1,\dots,5$). Hence
\[
I_{f,g}=(f^*\otimes g^*)(I_Z)=\big\langle \bar u_1,\dots,\bar u_5\big\rangle,\qquad
\bar u_i:=f^*(u_i)\otimes 1-1\otimes g^*(u_i).
\]

Consider the $5$–fold product $\bar u_1\smile\cdots\smile \bar u_5$. 
By expansion with the Koszul rule, its component in $H^5(Z)\otimes H^0(Z)$ is 
$(2\cdot 3\cdot 2\cdot 4\cdot 1)\,(u_1\smile\cdots\smile u_5)\otimes 1=48\,(u_1\smile\cdots\smile u_5)\otimes 1$,
and its component in $H^0(Z)\otimes H^5(Z)$ is 
$\pm(1\cdot 2\cdot 3\cdot 1\cdot 4)\,1\otimes (u_1\smile\cdots\smile u_5)=\pm 24\,1\otimes (u_1\smile\cdots\smile u_5)$.
Since both coefficients are nonzero in characteristic $0$, we have $\bar u_1\smile\cdots\smile \bar u_5\neq 0$.

On the other hand, any $6$–fold product of the $\bar u_i$’s vanishes: when expanding 
$\bar u_{i_1}\smile\cdots\smile \bar u_{i_6}$ each term is a tensor $(\alpha)\otimes(\beta)$ 
with $\alpha$ (resp.\ $\beta$) a cup product of some of the $u_j$’s; since there are only five 
degree–$1$ generators, either $\alpha$ or $\beta$ repeats an index and is $0$.
Therefore
$\mathrm{l.c.p.}\big(I_{f,g}\big)=5$ and $\mathrm{TC}(f,g)\geq 5.$
\end{example}

\section{The case of one map being a fibration: the Collaboration Principle}
A remarkable feature of the bivariate invariant is that 
the presence of a fibration on one side can simplify the synchronization problem. 
Intuitively, a fibration $g:Y\to Z$ provides a systematic way of lifting homotopies, so part of the motion planning task can be transferred to the lifting process. 
This effect is made precise in the following remarkable result.

\begin{theorem}[Collaboration Principle]\label{Colabor}
Let $f:X\to Z$ be a continuous map and let $g:Y\to Z$ be a surjective Hurewicz fibration. Then
$$
\mathrm{TC}(f,g) \leq \mathrm{TC}(f).
$$
In particular, whenever $\mathrm{TC}(f)$ is finite, so is  $\mathrm{TC}(f,g).$
\end{theorem}

\begin{proof}
Suppose an open subset $U\subseteq X\times Z$ together with a strict local section of $\pi _f:$
$$\xymatrix{
{U} \ar@{^(->}[rr]^{inc_{U}} \ar[dr]_{\sigma } & & {X\times Z} \\
 & {X^I} \ar[ur]_{\pi _f} & }$$
This means that $\sigma (x,z)(0)=x$ and $f(\sigma (x,z)(1))=z,$ for all $(x,z)\in U$. We define $V:=(id_X\times g)^{-1}(U)$ and 
$\sigma _1:V\to X^I$ as $\sigma _1(x,y):=\sigma (x,g(y)).$ Since $g$ is a fibration, we can take a lift $\widetilde{\Phi }$ in the following commutative diagram:
$$\xymatrix{
{V\times I\times \{1\}} \ar[rr]^{\varphi } \ar@{^(->}[d] & & {Y} \ar[d]^g \\
{V\times I\times I} \ar[rr]_{\Phi } \ar@{.>}[urr]^{\widetilde{\Phi }} & & {Z}  }$$
\noindent where $\varphi $ and $\Phi $ are defined as $\varphi (x,y,t,1):=y$ and $\Phi (x,y,t,s):=f(\sigma _1(x,y)(s)),$ respectively. We define
$$\sigma _2:V\to Y^I$$
\noindent as $\sigma _2(x,y)(s):=\widetilde{\Phi }(x,y,1,s)$. Then, it is straightforward to check that the pair $(\sigma_1(x,y), \sigma _2(x,y))$ is $(f,g)$-synchronized for all $(x,y)\in V$, that is, $f(\sigma _1(x,y)(s))=g(\sigma _2(x,y)(s))$, for all $s\in I$ and $(x,y)\in V.$ Moreover, $\sigma _1(x,y)(0)=x$ and $\sigma _2(x,y)(1)=y.$ In other words $\sigma '=(\sigma _1,\sigma _2):V\to \Delta _{f,g}^I$ defines a strict local section of $\pi _{f,g}:$
$$
\xymatrix{
{V} \ar@{^(->}[rr] \ar[dr]_{\sigma '} & & X \times Y \\
& \Delta_{f,g}^I \ar[ur]_{\pi_{f,g}} &
}
$$
Repeating this process to an open cover $\{U_i\}_{i=0}^n$ of $X\times Z$ that gives $\mathrm{TC}(f)=n$, we obtain an open cover $\{V_i\}_{i=0}^n$ of $X\times Y$ proving that $\mathrm{TC}(f,g)\leq n=\mathrm{TC}(f).$
\end{proof}

\begin{remark}
The Collaboration Principle can also be recovered as a particular case of 
Proposition~\ref{pre-comp-Pav}(2), which describes the behavior of bivariate topological 
complexity under precomposition by fibrations. In fact, if $g$ is a fibration, 
Proposition~\ref{pre-comp-Pav}(2) yields
\[
\mathrm{TC}(f,g)=\mathrm{TC}(f\circ \mathrm{id}_X,\mathrm{id}_Z\circ g)\leq \mathrm{TC}(f,\mathrm{id}_Z)=\mathrm{TC}(f).
\]
We have chosen to present a direct proof of Theorem~4.1, as it highlights the geometric 
mechanism underlying the collaboration phenomenon.
\end{remark}

Fix a fibration $g:Y\to Z$. In this setting, the bivariate topological complexity 
$\mathrm{TC}(f,g)$ defines an invariant of $f$ relative to $g$, which we denote by 
$\mathrm{TC}_g(f)$. Compared with the univariate invariant $\mathrm{TC}(f)$ of Pavešić, 
this relative version is technically simpler: its definition requires fewer constraints, 
and the presence of a fibration makes the construction of strict sections more transparent.

We now record an important rigidity property of the operator $\mathrm{TC}_g$ in the 
presence of a global section.

\begin{proposition}
Let $g:Y\to Z$ be a fibration admitting a global section. Then, for every map 
$f:X\to Z$, one has
\[
\mathrm{TC}_g(f)=\mathrm{TC}(f).
\]
\end{proposition}

\begin{proof}
The equality follows directly from Proposition~\ref{pre-comp-Pav}(3).
\end{proof}

\begin{remark}
This result shows that the relative invariant $\mathrm{TC}_g$ collapses to the classical
topological complexity of Pavešić whenever the synchronization constraint imposed by $g$
is trivial. In this sense, the bivariate framework detects only genuinely nontrivial
coordination phenomena.
\end{remark}

The following example shows that this collapse phenomenon does not hold in general when
the fibration $g$ does not admit a global section.

\begin{example}

Consider $f=\mathrm{pr}_1:S^2\times S^1 \to S^2$ the canonical projection and $g=p:PS^2 \to S^2$ the pointed path fibration. By using Proposition \ref{comparison}(3) and \cite[Cor. 4.16(1), Cor 4.9]{Sc}, we obtain $\mathrm{TC}(f)=\mathrm{TC}(S^2)=2$ and $\mathrm{TC}(g)=1$. Therefore, by Proposition \ref{bicolaboracion} we conclude that 
$$ \mathrm{TC}_g(f)=\mathrm{TC}(f,g)\leq 1 < 2 = \mathrm{TC}(f).$$
\end{example}

\bigskip
Now we prove that, whenever we have a fibration $g:Y\to Z$ we have that the operator $\mathrm{TC}_g$ is invariant under fibrewise homotopy equivalences. In order to do this we begin by analyzing how precomposition by a map affects the complexity of $(f,g)$. Compare with \cite[Theorem 3.5]{Pav}.

\begin{proposition}\label{precomp}
Let $X \stackrel{f}{\longrightarrow} Z \stackrel{g}{\longleftarrow} Y$ be a cospan of maps where $g$ is a fibration and consider a map $v:\widehat{X}\to X.$
\begin{enumerate}
\item If $v$ admits a homotopy section then $\mathrm{TC}_g(f)\leq \mathrm{TC}_g(f\circ v).$

\item If $v$ admits a homotopy section $u:X\to \widehat{X}$ such that $f\circ v\circ u=f,$ then $\mathrm{TC}_g(f\circ v)\leq \mathrm{TC}_g(f).$

\item If $v$ admits a homotopy section $u$ such that $f\circ v\circ u\simeq f$ and $f\circ v$ is a fibration, then $\mathrm{TC}_g(f\circ v)\leq \mathrm{TC}_g(f).$
    
\end{enumerate}
\end{proposition}

\begin{proof}${}$

(1) Take an open subset $U\subseteq \widehat{X}\times Y$ with a strict local section $\sigma =(\sigma _1,\sigma _2):U\to \Delta ^I _{f\circ v,g}$ of $\pi _{f\circ v,g}:$
$$\xymatrix{
{U} \ar@{^(->}[rr] \ar[dr]_{\sigma } & & {\widehat{X}\times Y} \\
 & {\Delta _{f\circ v,g}^I} \ar[ur]_{\pi _{f\circ v,g}} & }$$
 Now consider $u:X\to \widehat{X}$ a homotopy section of $v$ and take $H:X\times I\to X$ a homotopy with $H:v\circ u\simeq id_X.$ We define
 $V:=(u\times id_Y)^{-1}(U)\subseteq X\times Y$ and $\sigma '_1:V\to X^I$ as
 $$\sigma '_1(x,y)(t)=\begin{cases}
H(x,1-2t), & 0\leq t\leq \frac{1}{2} \\
v(\sigma _1(u(x),y)(2t-1)), & \frac{1}{2}\leq t\leq 1
 \end{cases}$$
Using the homotopy lifting property we take a lift in the following diagram
$$\xymatrix{
{V\times I\times \{1\}} \ar[rr]^{\varphi } \ar@{^(->}[d] & & {Y} \ar[d]^g \\
{V\times I\times I} \ar@{.>}[urr]^{\widetilde{\Phi }} \ar[rr]_{\Phi } & & {Z} }$$
\noindent where $\varphi $ and $\Phi $ are respectively defined as $\varphi (x,y,t,1):=y$ and $\Phi (x,y,t,s):=f(\sigma '_1(x,y)(s)).$ Defining
$\sigma '_2:V\to Y^I$ as $\sigma '_2(x,y)(s):=\widetilde{\Phi}(x,y,1,s)$ we obtain
$$\sigma '=(\sigma '_1,\sigma '_2):V\to \Delta ^I_{f,g}$$ a strict partial section of $\pi _{f,g}.$ 
Repeating this process to an open cover $\{U_i\}_{i=0}^n$ of $\widehat{X}\times Y$ that gives $\mathrm{TC}_g(f\circ v)=n$, we obtain an open cover $\{V_i\}_{i=0}^n$ of $X\times Y$ proving that $\mathrm{TC}_g(f)\leq n.$

(2) Similarly, let $U\subseteq X\times Y$ be an open subset together with a strict local section $\sigma =(\sigma _1,\sigma _2)$ of $\pi _{f,g}:$
$$\xymatrix{
{U} \ar@{^(->}[rr] \ar[dr]_{\sigma } & & {X\times Y} \\
 & {\Delta _{f,g}^I} \ar[ur]_{\pi _{f,g}} & }$$
We define $V:=(v\times id_Y)^{-1}(U)\subseteq \widehat{X}\times Y$ and $\tilde{\sigma }'_1:V\to X^I$ as
$\tilde{\sigma }'_1(\hat{x},y):=\sigma _1(v(\hat{x}),y).$ Now, if we explicit a homotopy $H:u\circ v\simeq id_{\widehat{X}}$, then we can define a map $\sigma '_1:V\to \widehat{X}^I$ as
$$
\sigma '_1(\hat{x},y)(t):= \begin{cases}
H(\hat{x},1-2t), & 0\leq t\leq \frac{1}{2} \\
u(\sigma _1(v(\hat{x}),y)(2t-1)), & \frac{1}{2}\leq t\leq 1\end{cases}$$ Again, using the homotopy lifting property of $g$ we take a lift
$$\xymatrix{
{V\times I\times \{1\}} \ar[rr]^{\varphi } \ar@{^(->}[d] & & {Y} \ar[d]^g \\
{V\times I\times I} \ar@{.>}[urr]^{\widetilde{\Phi }} \ar[rr]_{\Phi } & & {Z} }$$ \noindent where $\varphi (\hat{x},y,t,1):=y$ and $\Phi (\hat{x},y,t,s):=f(v(\sigma '_1(\hat{x},y)(s))).$ Then 
$$\sigma '_2:V\to Y^I$$ \noindent defined as $\sigma '_2(\hat{x},y)(s):=\widetilde{\Phi }(\hat{x},y,1,s)$, gives a strict local section $\sigma '=(\sigma '_1,\sigma '_2):V\to \Delta _{f\circ v,g}^I$ of $\pi _{f\circ v,g}.$ This gives $\mathrm{TC}_g(f\circ v)\leq \mathrm{TC}_g(f).$

(3) If $f\circ v$ is a fibration and $f\circ v\circ u\simeq f,$ then by the homotopy lifting property, there exists $u'\simeq u$ such that $f\circ v\circ u'=f.$ Since $v\circ u'\simeq id_X$ we can consider case (2) to obtain the result.
\end{proof}

Using similar techniques one can also show the following result:
\begin{proposition}
Let $X \stackrel{f}{\longrightarrow} Z \stackrel{g}{\longleftarrow} Y$ be a cospan of maps where $g$ is a fibration. If $v:\widehat{X}\to X$ is a fibration, then $\mathrm{TC}_g(f\circ v)\leq \mathrm{TC}_g(f).$  
\end{proposition}

\begin{corollary}
If $v:\widehat{X}\to X$ is a deformation retraction, then $\mathrm{TC}_g(f)=\mathrm{TC}_g(f\circ v),$ for any map $f:X\to Z$ and any fibration $g:X\to Z.$  
\end{corollary}

\begin{proof}
Let $u:X\to \widehat{X}$ be a map such that $v\circ u=id_X$ (so $u$ is, necessarily, an embedding) and $u\circ v\simeq id_{\widehat{X}}$. Then, by Proposition \ref{precomp} we have that $\mathrm{TC}_g(f)\leq \mathrm{TC}_g(f\circ v)$ and $\mathrm{TC}_g(f\circ v)\leq \mathrm{TC}_g(f).$
\end{proof}

\begin{corollary}
Let $X \stackrel{f}{\longrightarrow} Z \stackrel{g}{\longleftarrow} Y$ be a cospan of maps where $g$ is a fibration. If $f$ admits a homotopy section, then $\mathrm{TC}_g(f)\geq \mathrm{TC}(g).$ If, in addition, $f$ is a fibration, then $\mathrm{TC}_g(f)=\mathrm{TC}(g).$
\end{corollary}

\begin{proof}
By Proposition \ref{precomp} (1), we have
$\mathrm{TC}(g)=\mathrm{TC}(id_Z,g)\leq \mathrm{TC}(id_Z\circ f,g)=\mathrm{TC}_g(f).$

If $f$ is a fibration, then Theorem \ref{Colabor} (Collaboration Principle) holds and we have $\mathrm{TC}_g(f)=\mathrm{TC}(f,g)\leq \mathrm{TC}(g)$.
\end{proof}
Finally, we conclude this section with the following important result:  
\begin{theorem}
Let $g:Y\to Z$ be a fibration. Then the operator $\mathrm{TC}_g$ is invariant under fibrewise homotopy equivalences. In other words, if $f:X\to Z$ and $f':X'\to Z$ are fibrewise homotopy equivalent maps, we have
$$\mathrm{TC}_g(f)=\mathrm{TC}_g(f').$$
\end{theorem}

\begin{proof}
Consider $v:X\to X'$ a fibrewise homotopy equivalence between $f$ and $f'$:
$$\xymatrix{
{X} \ar[dr]_f \ar[rr]^v_{\simeq } & & {X'} \ar[ld]^{f'} \\
 & {Z} & }$$
 Then, by Proposition \ref{precomp} (1),
$$\mathrm{TC}_g(f)=\mathrm{TC}_g(f'\circ v)\geq \mathrm{TC}_g(f').$$
Similarly, taking a fibrewise homotopy inverse for $v$ we obtain  $\mathrm{TC}_g(f')\geq \mathrm{TC}_g(f)$ and we are done.
\end{proof}

\section{Scott-type bivariate topological complexity}\label{sec:Scott}

At this stage of the paper, we have studied a bivariate notion of topological complexity (of Pavešić type), which provides a measure of the minimal local information required to coordinate two systems whose dynamics are governed by a shared target space. However, although interesting in its own right, this notion is quite rigid and is not a homotopy invariant, as shown in Remark~\ref{ex-infty}. Consequently, it can be difficult to handle and compute.

In this section, we introduce a homotopy-invariant approach to bivariate topological complexity, called the Scott-type invariant, which is more flexible and computationally tractable. This construction provides a natural lower bound for the strict bivariate topological complexity. Its relevance stems from two main features: first, if both $f$ and $g$ are fibrations (even in the sense of Dold), then it coincides with the strict version; second, since it is defined in terms of homotopy sections rather than strict ones, it is more flexible and computationally accessible. We now proceed to give its definition, starting with the notion of $(f,g)$-sectional subset.

\begin{definition}
Given $X \stackrel{f}{\longrightarrow} Z \stackrel{g}{\longleftarrow} Y$ a cospan of continuous maps. and $A\subseteq X \times Y$, $A$ is called $(f,g)$-sectional if exists a map $\sigma: A\rightarrow Z^I$ such that $(f\times g) _{\vert_A} = \pi_Z \circ \sigma$, where $ \pi_Z $ is the usual path fibration in $Z$. So, given $(x,y)\in A$, it follows that 
$$(f\times g) _{\vert_A}(x,y)=(f(x),g(y))=(\pi_Z \circ \sigma)(x,y) =(\sigma(x,y)(0),\sigma(x,y)(1)). $$

This situation can be illustrated by the following diagram:
$$
\xymatrix{
{A} \ar[rr]^{(f\times g)_{\vert A}} \ar[dr]_{\sigma} & & {Z \times Z} \\ & Z^I   \ar[ur]_{\pi_Z} & } 
$$
In the above conditions, $\sigma$ is called a $(f,g)$-motion planner.
\end{definition}

The previous definition gives rise to the main notion of this section. 

\begin{definition}
Given $X \stackrel{f}{\longrightarrow} Z \stackrel{g}{\longleftarrow} Y$ continuous maps, the homotopy bivariate topological complexity of $f$ and $g$, denoted by $TC_H(f,g)$, is the least $k\in \mathbb{N}$ such that $X \times Y$ can be covered by $k+1$ $(f,g)$-sectional open subsets $U_0,\dots,U_k$. If no such $k$ exists, then we set $TC_H(f,g):=\infty .$
\end{definition}

\begin{remark}
$TC_H(f,g)$ will be also called Scott-type bivariate topological complexity. We will see that this is, indeed, a homotopy approach of  Pavešić type bivariate topological complexity. Such approach is, as we previously commented, easier to handle in practice, and often more amenable to explicit computation.

\end{remark}

We need first remember some facts about homotopy pullbacks. Recall from \cite{Mat} that given a homotopy commutative square
$$\xymatrix{
  P \ar[d]_{\alpha } \ar[r]^{\beta } & B \ar[d]^{g} \\
  A \ar[r]_{f} & C   } $$
\noindent equipped with a homotopy $H:f\circ \alpha \simeq g\circ \beta $,
another homotopy commutative square exists
$$\xymatrix{
  E_{f,g} \ar[d]_{p} \ar[r]^{q}  & B \ar[d]^{g} \\
A \ar[r]_{f} & C   }$$
\noindent where $E_{f,g}$ is the space defined as $E_{f,g}=\{(a,b,\lambda )
\in A\times B\times C^I \, ; \: \lambda (0)=f(a)\, , \:
\lambda (1)=g(b) \}$ and $p,q$ are the obvious restrictions of the
projections. This square is equipped with the homotopy $G$ given
by $(a,b,\lambda,t) \mapsto \lambda(t).$ Then, a map (referred to as the \textit{whisker map})
$w:P\rightarrow E_{f,g}$ exists, defined by
$w(x)=(\alpha (x), \beta (x), H(x,-))$ satisfying $p\circ w=\alpha
,$ $q\circ w=\beta $ and $G\circ (w\times id)=H.$ The first square is termed a
\textit{homotopy pullback} if the whisker map $w$ is a
homotopy equivalence; the second square is referred to as the \textit{standard homotopy pullback}
of $f$ and $g.$ It is well-known that, when $f$ or $g$ is a fibration, then their topological pullback (with the static homotopy) is a homotopy pullback. Moreover, the topological pullback $P_{f,g}$ is a strong deformation retract of $E_{f,g}.$

Returning to a cospan of continuous maps $X \stackrel{f}{\longrightarrow} Z \stackrel{g}{\longleftarrow} Y$, we can consider the standard homotopy pullback of $f$ and $g$,
$$
E_{f,g} = \{ (x, y, \lambda) \in X \times Y \times Z^I \mid \lambda(0) = f(x), \, \lambda(1) = g(y) \}.
$$
This space can be alternatively described via the following pullback diagram:
$$
\xymatrix{
E_{f,g} \ar[r]^{\overline{f\times g}} \ar[d]_{\theta _{f,g}=\overline{\pi}_Z} & Z^I \ar[d]^{\pi_Z} \\
X \times Y \ar[r]_{f \times g} & Z \times Z.
}
$$

Observe that, in order to emphasize the role of the maps $f$ and $g,$ we will rather denote by $\theta _{f,g}$ the induced map $\overline{\pi }_Z$ of $\pi _Z$ in this pullback. We can establish the following result:

\begin{theorem}
Given $X \stackrel{f}{\longrightarrow} Z \stackrel{g}{\longleftarrow} Y$ a cospan of continuous maps, the following are equivalent for a subset $A\subseteq X\times Y$:

\begin{enumerate}
    \item There is a local section $s:A\longrightarrow E_{f,g}$ of $\theta _{f,g}.$

    \item There is an $(f,g)$-motion planner $\sigma:A\longrightarrow Z^I$.

    \item The map $(f\times g)_{\vert A} $ can be deformed into $\Delta_Z(Z)$, where $\Delta _Z:Z\to Z\times Z$ denotes the diagonal map in $Z.$

    \item There is a map $h:A\longrightarrow Y$ such that 
    $ (f\times g)_{\vert A} $ is homotopic to $ (g\times g) \circ \Delta_{Y} \circ h$.
    
    \item There is a map $\widehat{h}:A\longrightarrow X$ such that $ (f\times g)_{\vert A} $ is homotopic to $ (f\times f) \circ \Delta_{X} \circ \widehat{h}$.
\end{enumerate}
\end{theorem}

\begin{proof}
Let us first check that (1) implies (2).  
Let $s:A\longrightarrow E_{f,g}$ be the existing local section of $\theta _{f,g}$; thus $\theta _{f,g} \circ s=inc_A$. We define $\sigma:A \longrightarrow Z^I$ given by $\sigma :=(\overline{f\times g}) \circ s $. Then, $\sigma$ is an $(f,g)$-motion planner. Indeed, taking into account the commutativity of the pullback that defines $E_{f,g}$, it follows that:
$$\pi_Z \circ \sigma =\pi_Z\circ (\overline{f\times g}) \circ s=(f\times g)\circ \theta _{f,g} \circ s=(f\times g)\circ inc_A=(f\times g)_{\vert_A}.$$

\medskip
For (2) implies (1), 
let $\sigma:  A \longrightarrow Z^I$ be an $(f,g)$-motion planer. By the universal property of the pullback 
$$
\xymatrix{
A \ar@/^1pc/[drr]^{\sigma} \ar@/_1pc/[ddr]_{inc_A} \ar@{.>}[dr]^{s} & & \\
 & E_{f,g} \ar[r]^{\overline{f\times g}} \ar[d]_{\theta _{f,g}} & Z^I \ar[d]^{\pi_Z} \\
 & X \times Y \ar[r]_{f\times g} & Z \times Z.
}
$$
\noindent there exists a map $s:A\longrightarrow E_{f,g}$ satisfying $(\overline{f\times g})\circ s=\sigma $ and $\theta _{f,g}\circ s=inc_A $ obtaining this way that $s$ is a local section of $\theta _{f,g}$. 

\medskip
Now we check (2) implies (3):
Let $\sigma:  A \longrightarrow Z^I$ be an $(f,g)$-motion planer. We define the homotopy $H:A\times I\longrightarrow Z\times Z$ given by $H(x,y,t):= (\sigma(x,y)(t),g(y))$. Then $H(x,y,0)=(\sigma(x,y)(0),g(y))=(f(x),g(y))=(f\times g)_{\vert A}(x,y)$ and $H(x,y,1)=(\sigma(x,y)(1),g(y))=(g(y),g(y))\in \Delta _Z(Z).$ 

\medskip
In order to prove (3) implies (2) consider a homotopy $H:A\times I\longrightarrow Z\times Z$ such that $H(x,y,0)=(f(x),g(y))$ and $H(x,y,1)\in \Delta _Z(Z)$, we define $\sigma: A\longrightarrow Z^I$ as follows:
$$ \sigma(x,y)(t) = \begin{cases}
pr_1(H(x,y,2t)), & \text{if } 0\leq t \leq \frac{1}{2}\\
pr_2(H(x,y,2-2t)),  & \text{if } \frac{1}{2} \leq t \leq 1 \\
\end{cases}$$
It is straightforward to check that $\sigma$ is an $(f,g)$-motion planner.

\medskip
For (2) implies (4), let us take the homotopy $H$ defined in (2)$\Rightarrow$ (3); there we checked that $H(x,y,1)=(g(y),g(y))$. We observe that $$(g(y),g(y)=(g\times g)(y,y)=((g\times g)\circ \Delta_Y \circ pr_2 )(x,y),$$ \noindent where $pr_2: X\times Y \longrightarrow Y$ is the canonical projection and $\Delta_Y: Y \longrightarrow Y\times Y$ the diagonal map on $Y$. Therefore, taking $h=pr_2$, we conclude that $H$ is a homotopy between $(f\times g)_{\vert A}$ and the composition $(g\times g)\circ \Delta_Y \circ pr_2.$

\medskip Analogously to the previous reasoning, in order to prove that (2) implies (5), we can take the homotopy $G(x,y,t):= (\sigma(x,y)(t),f(x))$ between $(f\times g)_{\vert A}$ and $(f\times f)\circ \Delta_X \circ pr_1.$

\medskip
Let us show that (4) implies (3):
Consider a map $h:A\longrightarrow Y$ satisfying that $ (f\times g)_{\vert A} $ is homotopic to $ (g\times g) \circ \Delta_{Y} \circ h$.  Then $(g\circ \Delta_{Y}\circ h )(x,y)=(g(h(x,y)),g(h(x,y)))\in \Delta _Z(Z).$ Hence, $ (g\times g) \circ \Delta_{Y} \circ h$ can be deformed into $\Delta_{Z}(Z)$.  

\medskip
Finally, the proof of (5) implies (3) is analogous to (4) implies (3).

\end{proof}

The map $\theta _{f,g}=\overline{\pi }_Z$ is a fibration since it is induced by the path fibration $\Pi_Z$ in the above pullback diagram. It allows us to establish the following key result

\begin{corollary}\label{pullback-def}
Given two continuous maps $f:X \rightarrow Z$ and $g:Y \rightarrow Z$, then $$TC_H(f,g)=secat(\theta _{f,g})$$ 
\end{corollary}

\begin{proof}
Just consider the equivalence between (1) and (2) of the previous theorem.
\end{proof}

As a consequence, when \(X\), \(Y\), and \(Z\) are ANR spaces, one can define a generalized version of \(\mathrm{TC}_H\) by allowing arbitrary (not necessarily open) covers:

\begin{proposition}
Let \(X \xrightarrow{f} Z \xleftarrow{g} Y\) be a cospan of ANR spaces. Then \(\mathrm{TC}_H(f,g)\) is the least integer \(k \in \mathbb{N}\) (or infinity) such that \(X \times Y\) admits a cover by \(k+1\) \((f,g)\)-sectional subsets \(U_0,\dots,U_k\), not necessarily open.
\end{proposition}

\begin{proof}
This follows directly from \cite[Proposition 5.10]{GC}, together with Corollary~\ref{pullback-def} above, which shows that \(\mathrm{TC}_H(f,g)\) coincides with the relative sectional category \(\mathrm{secat}_{f\times g}(\pi_Z)\) in the sense of J.~González, M.~Grant, and L.~Vandembroucq \cite{G-G-V}.
\end{proof}

One can also see that $\mathrm{TC}_H$ is symmetric, since the maps 
$f \times g$ and $g \times f$ are related by the homeomorphism 
$X \times Y \cong Y \times X$ that swaps the factors. Moreover, it is also immediate to check that $\mathrm{TC}_H(id_X,id_X)=\mathrm{TC}(X),$ for every space $X.$

\medskip
Unlike Pave\v{s}i\'c-type bivariate topological complexity, this approach is much more robust, as the following result displays:

\begin{corollary}
Let $f,f'\colon X \to Z$ and $g,g'\colon Y \to Z$ be continuous maps such that $f\simeq f'$ and $g\simeq g'$. Then
\[
\mathrm{TC}_H(f,g)=\mathrm{TC}_H(f',g').
\]
\end{corollary}

\begin{proof}
Since $f\simeq f'$ and $g\simeq g'$, it follows that $f\times g\simeq f'\times g'$.  
As $\pi_Z \colon Z^I\to Z\times Z$ is a fibration, the pullbacks of $\pi_Z$ along $f\times g$ and $f'\times g'$ are fibrewise homotopy equivalent (see, for instance, \cite[Chapter~2, §8, Theorem~14]{Sp}).  
Since the sectional category is invariant under (fibrewise) homotopy equivalence, we conclude that
$\mathrm{TC}_H(f,g)=\mathrm{TC}_H(f',g'),$
as required.
\end{proof}

\section{Comparison of $\mathrm{TC}_H(f,g)$ and $\mathrm{TC}(f,g)$}

In the previous sections we introduced two different formulations of bivariate 
topological complexity: the \emph{strict version} $\mathrm{TC}(f,g)$, defined via genuine 
sections of the map $\pi_{f,g}$, and the \emph{homotopy-invariant version} 
$\mathrm{TC}_H(f,g)$, defined through homotopy sections of the fibration $\theta_{f,g}$. 
Although both invariants measure the cost of coordinating the maps $f$ and $g$, 
they behave quite differently: $\mathrm{TC}(f,g)$ is often rigid and may take the 
value~$\infty$, while $\mathrm{TC}_H(f,g)$ is more flexible but may lose information. 
It is therefore natural to ask how these two quantities are related. 
In this section we establish comparison results between them, showing in particular 
that they coincide under favourable hypotheses (notably when $f$ and $g$ are fibrations), 
but can differ drastically in general.

\begin{lemma}\label{fib}
If $f$ and $g$ are fibrations, then $\pi _{f,g}$ is also a fibration, and therefore $$\mathrm{TC}(f,g)=\mathrm{secat}(\pi _{f,g}).$$
\end{lemma}

\begin{proof}
Observe that if $f : X \to Z$ and $g : Y \to Z$ are fibrations, then the pullback diagram
$$\xymatrix{
\Delta _{f,g} \ar[d]_{\overline{g}} \ar[r]^{\overline{f}} & Y \ar[d]^g \\
X \ar[r]_f & Z} $$
induces pullback maps $\overline{f}$ and $\overline{g}$, which are themselves fibrations. As a result, their product map $\overline{g} \times \overline{f} : \Delta _{f,g} \times \Delta _{f,g} \to X \times Y$ is also a fibration. Now consider the canonical path fibration $\pi : \Delta _{f,g}^I \to \Delta _{f,g} \times \Delta _{f,g},$
which assigns to a path $\gamma = (\alpha, \beta)$ its endpoints $(\gamma(0), \gamma(1))$. The composite $\pi_{f,g} = (\overline{g} \times \overline{f}) \circ \pi : \Delta _{f,g}^I \to X \times Y$
is then a fibration as well. This situation is summarized in the following commutative diagram:
$$
\xymatrix{
\Delta _{f,g}^I \ar@{.>}[rr]^{\pi_{f,g}} \ar[dr]_{\pi} & & X \times Y \\
& \Delta _{f,g} \times \Delta _{f,g} \ar[ur]_{\overline{g} \times \overline{f}} &
}$$    
\end{proof}

Recall that a \emph{Dold fibration} (or weak fibration) is a map  $p:E \to B$ satisfying the weak covering homotopy property (see \cite{D} for details). Such maps are completely characterized by being fibrewise homotopy equivalent to Hurewicz fibrations. 

One has the following interesting result:

\begin{lemma}
If $p$ is a Dold fibration, then $\mathrm{sec}(p) = \mathrm{secat}(p)$.    
\end{lemma}

\begin{proof}
Just take into account Proposition \ref{triangle}(2), Proposition \ref{invariance}(1) and the fact that $\mathrm{sec}$ and $\mathrm{secat}$ coincide on the class of Hurewicz fibrations.     
\end{proof}

In addition, Dold fibrations are stable under pullbacks, products and compositions. Consequently, Lemma~\ref{fib} remains valid when both $f$ and $g$ are Dold fibrations. One of the main motivations for considering Dold fibrations in the study of 
topological complexity is their greater flexibility compared to Hurewicz fibrations. In particular, the class of Dold fibrations is broad enough to include certain continuous maps that fail to be Hurewicz fibrations, 
typically due to the presence of singularities or degeneracies, yet still exhibit good fibrewise homotopical behaviour. This makes it possible to 
extend the framework of sectional category and topological complexity beyond the more restrictive setting of Hurewicz fibrations, in line with the perspective advocated by Pave\v{s}i\'c.

\medskip

On the other hand, if either $f$ or $g$ is a fibration, then the ordinary pullback $\Delta _{f,g}$ is a homotopy pullback. Consequently, the associated whisker map $w : \Delta _{f,g} \to E_{f,g}$ is a homotopy equivalence. In fact, the following diagram commutes:
$$
\xymatrix{
\Delta _{f,g} \ar[rr]^w_{\simeq} \ar@{^(->}[dr]_{\mathrm{inc}} & & E_{f,g} \ar[dl]^{\theta _{f,g}} \\
& X \times Y . &
}
$$

\begin{proposition}\label{comparison}
Let $X \stackrel{f}{\longrightarrow} Z \stackrel{g}{\longleftarrow} Y$ be any cospan of continuous maps. 
Then
\begin{enumerate}
\item $\mathrm{TC}_H(f,g) \;\leq\; \mathrm{TC}(f,g);$
\item If at least one of $f$ or $g$ is a fibration, then
$\mathrm{TC}_H(f,g) \;=\; \mathrm{secat}(\pi_{f,g});
$

\item If both $f$ and $g$ are fibrations, then
$$\mathrm{TC}_H(f,g) \;=\; \mathrm{TC}(f,g).$$
\end{enumerate}

\end{proposition}

\begin{proof}
In order to prove (1) consider the following commutative diagram
$$\xymatrix{
{\Delta _{f,g}^I} \ar[d]_{\pi _{f,g}} & {\Delta _{f,g}} \ar[d]^{inc} \ar[l]^{\simeq }_{c} \ar@{^(->}[d] \ar[r]^{w} & {E_{f,g}} \ar[d]^{\theta _{f,g}} \\
{X\times Y} \ar@{=}[r] & {X\times Y} \ar@{=}[r] & {X\times Y.}
}$$
Here, $w$ is the whisker map and $c$ is the map that assigns to each point $(x,y)\in \Delta _{f,g}$ the constant path at $(x,y),$ $c_{(x,y)}$, which is a homotopy equivalence. We therefore obtain a chain of equalities:
$$\mathrm{TC}_H(f,g)=\mathrm{secat}(\theta _{f,g})\leq \mathrm{secat}(inc)=\mathrm{secat}(\pi _{f,g})\leq \mathrm{sec}(\pi _{f,g})=\mathrm{TC}(f,g).$$

Now we check (2). If one of the maps is a fibration, then $w$ is a homotopy equivalence, and one can directly see that $\mathrm{TC}_H(f,g)=\mathrm{secat}(\pi_{f,g}).$

Finally we prove (3). In this case the whisker map is also a homotopy equivalence and, by Lemma \ref{fib}, $$\mathrm{TC}_H(f,g)=\mathrm{secat}(\theta _{f,g})=\mathrm{secat}(inc)=\mathrm{secat}(\pi _{f,g})=\mathrm{TC}(f,g).$$

\end{proof}

\begin{remark}
Again, Proposition \ref{comparison} above also holds if we consider Dold fibrations instead of fibrations. We leave the details to the reader.
\end{remark}

\begin{example}
The invariants $\mathrm{TC}(f,g)$ and $\mathrm{TC}_H(f,g)$ may differ, even drastically. Now we give a list of different situations that illustrates this behavior:
\begin{enumerate}
\item 
Let $f=c_{z}:X\to Z$ and $g=c_{z'}:Y\to Z$ be constant maps at distinct points 
$z\neq z'$ belonging to the same path-component of $Z$. 
As observed in Remark \ref{ex-infty}, in this case the pullback defining the strict invariant is empty, hence 
$\mathrm{TC}(f,g)=\infty$. 
In contrast, the choice of a single path $\gamma:[0,1]\to Z$ with $\gamma(0)=z$ and $\gamma(1)=z'$ 
yields a continuous global $(f,g)$-motion planner 
$\sigma :X\times Y\to Z^I$, showing that $\mathrm{TC}_H(f,g)=0$. 
This example highlights the rigidity of the strict version versus the flexibility of the homotopy-invariant formulation.

\item A different iconic phenomenon appears in Example \ref{iconic-example}. 
There, for the maps $f,g:S^1\to \mathbb{C}$ defined by $f(z)=z^2$ and $g(z)=-z^2$, we established that $\mathrm{TC}(f,g)=1$. 
However, since $\mathbb{C}$ is convex, one can choose a linear homotopy between $f(z)$ and $g(w)$ for every $(z,w)\in S^1\times S^1$. 
Explicitly, a global $(f,g)$-motion planner is given by
\[
\sigma(z,w)(t) = (1-t)f(z) + t g(w) = (1-t)z^2 -tw^2,
\]
which shows that $\mathrm{TC}_H(f,g)=0$. 

\item Finally, the case of inclusions $\mathrm{inc}_1:A_1\hookrightarrow X$, $\mathrm{inc}_2:A_2\hookrightarrow X$  provides a ``bridge'' example. 
Here the strict invariant recovers the classical topological complexity,
\[
\mathrm{TC}(\mathrm{inc}_1,\mathrm{inc}_2)=\begin{cases}
\infty, & \mbox{if}\hspace{4pt} A_1\neq A_2 \\
\mathrm{TC}(A), & \mbox{if}\hspace{4pt} A_1=A_2=A
\end{cases}
\]
while the homotopic version yields the relative complexity in the sense of Farber \cite{Far2} of $A_1\times A_2$ in $X:$
\[
\mathrm{TC}_H(\mathrm{inc}_1,\mathrm{inc}_2)=\mathrm{TC}_X(A_1\times A_2).
\]
For instance, with $\mathrm{inc}_1=\mathrm{inc}_2=inc:S^2\hookrightarrow \mathbb{R}^3$, the strict invariant is $\mathrm{TC}(S^2)=2$  while the homotopic one vanishes, $\mathrm{TC}_H(\mathrm{inc},\mathrm{inc})=0$, because $\mathbb{R}^3$ is convex. 
Thus both numbers are finite but encode completely different information: the intrinsic complexity of $A$ versus the ambient cost inside $X$.
\end{enumerate}
Together, these examples show that $\mathrm{TC}(f,g)$ and $\mathrm{TC}_H(f,g)$ are genuinely different invariants, complementary in spirit and in applications.
\end{example}

\section{Relationship between $\mathrm{TC}_H$ and the homotopic distance}
Another noteworthy feature of this homotopy invariant is that it can be expressed
in terms of the homotopic distance of Macías-Virgós and Mosquera-Lois \cite{M-M}.
Given a cospan of continuous maps \(X \xrightarrow{f} Z \xleftarrow{g} Y\), let
\(\mathrm{pr}_1:X\times Y\to X\) and \(\mathrm{pr}_2:X\times Y\to Y\) denote the canonical projections.
Since \(f\times g=(f\circ \mathrm{pr}_1,\,g\circ \mathrm{pr}_2):X\times Y\to Z\times Z\), it follows
from \cite[Theorem 2.7]{M-M} together with Corollary~\ref{pullback-def} that
\[
\TC_H(f,g)=\mathrm{D}\big(f\circ \mathrm{pr}_1,\,g\circ \mathrm{pr}_2\big).
\]

This way, many properties coming from the homotopic distance can be transfered to $\mathrm{TC}_H.$ We do not want to exhibit all these properties here, but we can give an example which provides the subadditivity of the product for homotopy bivariate topological complexity:

\begin{proposition}
Let $X \xrightarrow{f} Z \xleftarrow{g} Y$  and  $X' \xrightarrow{f} Z' \xleftarrow{g} Y'$ be cospans of continuous maps such that the product $X\times X'\times Y\times Y'$ is normal. Then we have
$$\mathrm{TC}_H(f\times f',g\times g')\leq \mathrm{TC}_H(f,g)+\mathrm{TC}_H(f',g').$$
\end{proposition}

\begin{proof}
Just consider \cite[Theorem 3.19]{M-M} and the formula $\TC_H(f,g)=\mathrm{D}\big(f\circ \mathrm{pr}_1,\,g\circ \mathrm{pr}_2\big)$ for any cospan $(f,g)$ of maps. We leave the details to the reader.   
\end{proof}

\bigskip
When restricting to maps \(X\to Z\) (that is, with common domain \(X\)), we obtain the
following comparison between the homotopic distance and the homotopy bivariate
topological complexity. More precisely, we prove that, under fairly mild conditions, $\TC_H(f,g)$ takes values in the closed interval $[\mathrm{D}(f,g),\mathrm{D}(f,g)+\TC(X)],$ which has length $\mathrm{TC}(X).$ Thus, $\TC(X)$ provides an upper bound for the possible deviation of $\TC_H(f,g)$ from $\mathrm{D}(f,g):$ 

\begin{proposition}\label{rel-dist}
Let \(f,g:X\to Z\) be maps. Then
\[
\mathrm{D}(f,g)\leq \TC_H(f,g).
\]
If, in addition, \(X\times X\) is normal (e.g.\ \(X\) is metrizable), then
\[
\mathrm{D}(f,g)\leq \TC_H(f,g)\leq \mathrm{D}(f,g)+\TC(X).
\]
\end{proposition}

\begin{proof}
Consider the consecutive pullback squares
\[
\xymatrix{
\mathcal{P}(f,g) \ar[r] \ar[d]_{p_{f,g}} & E_{f,g}  \ar[r] \ar[d]^{\theta _{f,g}} & {Z^I} \ar[d]^{\pi _Z} \\
X \ar[r]_{\Delta_X} & X\times X \ar[r]_{f\times g} & {Z\times Z}
}
\]
Since the composition of pullbacks is a pullback and $(f\times g)\circ \Delta _X=(f,g),$ then by \cite[Theorem 2.7]{M-M} we have \(\mathrm{D}(f,g)=\mathrm{secat}(p_{f,g})\).
Since sectional category is non-increasing under pullback (see Proposition \ref{pullback}), we obtain
\[
\mathrm{D}(f,g)=\mathrm{secat}(p_{f,g})\le \mathrm{secat}(\theta _{f,g})=\TC_H(f,g).
\]

For the proof of the second statement, use the triangle inequality in \cite[Prop. 3.16]{M-M}, together with \cite[Prop. 2.6]{M-M} and the behavior of the homotopic distance under pre- and post-composition \cite[Prop. 3.3 and Prop. 3.1]{M-M}:

\[
\begin{aligned}
\TC_H(f,g)
&= \mathrm{D}\big(f\circ \mathrm{pr}_1,\,g\circ \mathrm{pr}_2\big) \\
&\le \mathrm{D}\big(f\circ \mathrm{pr}_1,\,f\circ \mathrm{pr}_2\big)
   + \mathrm{D}\big(f\circ \mathrm{pr}_2,\,g\circ \mathrm{pr}_2\big) \\
&\le \mathrm{D}(\mathrm{pr}_1,\mathrm{pr}_2)+\mathrm{D}(f,g)
 \;=\; \TC(X)+\mathrm{D}(f,g).
\end{aligned}
\]
\end{proof}

\begin{corollary}
Let $f,g:X\to Z$ be continuous maps where $X\times X$ is normal. If $f\simeq g,$ then 
$$\mathrm{TC}_H(f,g)\leq \mathrm{TC}(X).$$    
\end{corollary}

\begin{remark}
Observe that in the corollary above we recover the result $\mathrm{TC}_{Sc}(f)\leq \mathrm{TC}(X)$, where $\mathrm{TC}_{Sc}(f):=\mathrm{TC}_H(f,f)$ denotes the topological complexity of $f$ in the sense of Scott \cite{Sc}.
\end{remark}

As an immediate consequence, when we restrict to maps \(X\to Z\) we obtain a triangle
inequality for \(\TC_H\).

\begin{proposition}
Let \(f,g,h:X\to Z\) be maps, with \(X\times X\) normal (e.g.\ \(X\) metrizable). Then
\[
\TC_H(f,g)\leq \TC_H(f,h)+\TC_H(h,g).
\]
\end{proposition}

\begin{proof}
By the triangle inequality for the homotopic distance (see \cite[Prop. 3.16]{M-M}),
\[
\TC_H(f,g)=\mathrm{D}\big(f\circ \mathrm{pr}_1,\,g\circ \mathrm{pr}_2\big)
\leq \mathrm{D}\big(f\circ \mathrm{pr}_1,\,h\circ \mathrm{pr}_1\big)
+\mathrm{D}\big(h\circ \mathrm{pr}_1,\,g\circ \mathrm{pr}_2\big).
\]
By \cite[Prop. 3.3]{M-M} and Proposition~\ref{rel-dist},
\[
\mathrm{D}\big(f\circ \mathrm{pr}_1,\,h\circ \mathrm{pr}_1\big)\leq \mathrm{D}(f,h)\leq \TC_H(f,h),
\]
and \(\mathrm{D}\big(h\circ \mathrm{pr}_1,\,g\circ \mathrm{pr}_2\big)=\TC_H(h,g)\), which proves the claim.
\end{proof}

It should be emphasized, however, that \(\TC_H(f,g)\) is not a distance
on the set of maps \(X\to Z\) with fixed domain: in general \(\TC_H(f,f)\neq 0\).
Indeed, for \(f=\mathrm{id}_X\) one has \(\TC_H(\mathrm{id}_X,\mathrm{id}_X)=\TC(X)\), which
vanishes only when \(X\) is contractible. In fact, as we have previously commented, \(\TC_H(f,f)\) agrees precisely with
Scott’s topological complexity of a map \cite{Sc}. 

For a map $f:X\to Z$, we consider the notation 
$$\mathrm{TC}_H(f):=\mathrm{TC}_H(f,id_Z).$$
Note that $\mathrm{TC}_H(f)$ agrees with the so called \emph{mixed topological complexity} of $f$ in the sense of Scott \cite{Sc}. In view of Proposition \ref{comparison} it can also be defined as
$$\mathrm{TC}_H(f)=\mathrm{secat}(\pi _f)$$ \noindent where $\pi _f:X^I\to X\times Z$ is Pavesic's map, which defines $\mathrm{TC}(f).$
Moreover, $\mathrm{TC}_H(f)\leq \mathrm{TC}(f).$ We have the following corollary:

\begin{corollary}
Let $f,g:X\to Z$ two maps with the same domain $X$ such that $X\times X$ is normal. Then
$$\mathrm{TC}_H(f,g)\leq \mathrm{TC}_H(f)+\mathrm{TC}_H(g).$$
\end{corollary}

\begin{proof}
Just consider $h=id_Z$ and the triangle inequality.    
\end{proof}

\medskip

\section{Bounds for $\mathrm{TC}_H$}

Let $X \xrightarrow{f} Z \xleftarrow{g} Y$ be a cospan, and set 
$F = f\circ \mathrm{pr}_1$ and $G = g\circ \mathrm{pr}_2 : X\times Y \to Z$, 
where $\mathrm{pr}_1,\mathrm{pr}_2$ are the projections from $X\times Y$ onto the first and second factor, respectively. 
Recall that $\mathrm{TC}_H(f,g) = D(F,G)$, where $D$ denotes the homotopic distance. 
Following Macías-Virgós and Mosquera-Lois, define
\[ 
\mathcal{J}(F,G) \;=\; \mathrm{im}\!\big(F^* - G^*: H^*(Z;R) \to H^*(X\times Y;R)\big) 
\;\subseteq\; H^*(X\times Y;R).
\] 
Their cohomological theorem \cite[Theorem 5.2]{M-M} asserts that
\[
\mathrm{l.c.p.}\,\mathcal{J}(F,G) \;\leq\; D(F,G) \;=\; \mathrm{TC}_H(f,g).
\]

On the other hand, the strict bivariant topological complexity $\mathrm{TC}(f,g)$ 
admits a cohomological bound in terms of the submodule
\[
I_{f,g} \;=\; (f\times g)^*\!\big(\ker \Delta^*_Z\big) \;\subseteq\; H^*(X\times Y;R),
\]
where $\Delta _Z: Z \to Z\times Z$ is the diagonal and 
$\pi_1,\pi_2: Z\times Z \to Z$ are the coordinate projections. In particular,
\[
\mathrm{l.c.p.}(I_{f,g}) \;\leq\; \mathrm{TC}(f,g).
\]

For every $\alpha\in H^*(Z;R)$ we have 
\[
(F,G)^*(\pi_1^*\alpha - \pi_2^*\alpha) \;=\; F^*\alpha - G^*\alpha \;\in\; \mathcal{J}(F,G).
\]
Moreover, $\ker\Delta_Z^*$ is generated (as an $R$–submodule of $H^*(Z\times Z;R)$) by finite sums of finite products of classes of the form $\pi_1^*\alpha-\pi_2^*\alpha$. Hence, $I_{f,g} = (f\times g)^*(\ker\Delta_Z^*)$ is the $R$–submodule generated by all finite products of elements of $\mathcal{J}(F,G).$
Since the lower cup–length depends only on the maximal number of factors in a nontrivial product, it follows that
\[
\mathrm{l.c.p.}\,\mathcal{J}(F,G)\;=\;\mathrm{l.c.p.}(I_{f,g}).
\]

Consequently, the cohomological lower bounds arising from the homotopic and the 
strict versions always coincide:
\[
\mathrm{l.c.p.}\,\mathcal{J}(F,G)\;=\;\mathrm{l.c.p.}(I_{f,g})
\;\leq\;\mathrm{TC}_H(f,g)\;\leq\;\mathrm{TC}(f,g).
\]

 In analogy with the strict Pavešić version, the homotopic setting also satisfies a \emph{Collaboration Principle}, but here without any fibrancy assumptions.  

\begin{proposition}\label{bicolaboracion}
Let $X \xrightarrow{f} Z \xleftarrow{g} Y$ be a cospan of maps. Then
$$
\mathrm{TC}_H(f,g)\le \min\{\mathrm{TC}_H(f),\,\mathrm{TC}_H(g)\}.
$$
\end{proposition}

\begin{proof}
Since $(f\times \mathrm{id}_Z)\circ (\mathrm{id}_X\times g)=f\times g$, taking consecutive pullbacks we have
$$\xymatrix{
{E_{f,g}} \ar[d]_{\theta _{f,g}} \ar[r] & {E_{f,\mathrm{id}_Z}} \ar[d]_{\theta _{f,\mathrm{id}_Z}} \ar[r] & {Z^I} \ar[d]^{\pi _Z} \\
{X\times Y} \ar[r]_{\mathrm{id}_X\times g} & {X\times Z} \ar[r]_{f\times \mathrm{id}_Z} & {Z\times Z}
}$$
Therefore, by Proposition\ref{pullback}
$$
\mathrm{TC}_H(f,g)=\mathrm{secat}(\theta _{f,g})\leq \mathrm{secat}(\theta _{f,\mathrm{id}_Z})
=\mathrm{TC}_H(f,\mathrm{id}_Z)=\mathrm{TC}_H(f).
$$
Similarly, $\mathrm{TC}_H(f,g)\leq \mathrm{TC}_H(g).$
\end{proof}

\begin{remark}
In the strict Pavešić version an analogous inequality requires a fibrancy assumption on one leg (e.g. $g$ a fibration), yielding 
$\mathrm{TC}(f,g)\le \mathrm{TC}(f).$
The homotopic variant is thus more flexible.
\end{remark}

\bigskip

Before turning to finer bounds, we anchor the homotopy bivariate invariant between familiar LS–type quantities. On the one hand, $\TC_H(f,g)$ cannot be smaller than the category of the individual maps $f$ and $g$, so it faithfully reflects the intrinsic planning cost already present in each factor. On the other hand, $\TC_H(f,g)$ is bounded above by the category of the combined map $f\times g$ and by the ordinary topological complexity of the target $Z$. These estimates situate $\TC_H(f,g)$ within a well–understood range and will serve repeatedly as a baseline for computations and comparisons.

\begin{proposition}\label{gen-bounds}
Let $X \xrightarrow{f} Z \xleftarrow{g} Y$ be a cospan of maps. Then
the following inequalities hold:
\begin{enumerate}
\item $\mathrm{TC}_H(f,g)\leq \mathrm{TC}(Z).$

\item If $Z$ is path-connected, then
$$ \max\{\mathrm{cat}(f),\mathrm{cat}(g) \}\leq \mathrm{TC}_H(f,g)\leq \mathrm{cat}(f\times g).$$
\end{enumerate}
\end{proposition}
\begin{proof}
(1) By Proposition\ref{pullback}, we clearly have 
$$\mathrm{TC}_H(f,g)=\mathrm{secat}(\theta _{f,g})\leq \mathrm{secat}(\pi _Z)=\mathrm{TC}(Z).$$
(2) In order to prove that $\mathrm{cat}(f)\leq \mathrm{TC}(f,g)$, suppose $\mathrm{TC}_H(f,g)=n$ and consider $\{U_i\}_{i=1}^n$ an open cover of $X\times Y$ with $(f,g)$-motion planners $\sigma _i:U_i\to Z^I$. If we fix a point $y_0\in Y$, we define $V_i:=\{x\in X:(x,y_0)\in U_i\},$ which is an open subset of $X$, and a homotopy
$$H_i:V_i\times I\to Z$$ \noindent given as $H_i(x,t):=\sigma _i(x,y_0)(t).$ We then have $H_i(x,0)=f(x)$ and $H_i(x,1)=g(y_0),$ proving that $f_{|V_i}$ is nullhomotopic. Therefore, $\mathrm{cat}(f)\leq \mathrm{TC}_H(f,g).$ The inequality $\mathrm{cat}(g)\leq \mathrm{TC}(f,g)$ is similarly proved by symmetry.

Now, suppose that $\mathrm{cat}(f\times g)=n$ and consider an open cover $\{V_j\}_{j=0}^n$ of $X\times Y$ such that every restriction $(f\times g)_{|V_j}$ is nullhomotopic. Since $Z$ is path-connected, for a fixed element $z_0\in Z$ we can consider homotopies
$$F_j:V_j\times I\to Z\times Z$$
\noindent such that $F_j(x,y,0)=(f(x),g(y))$ and $F_j(x,y,1)=(z_0,z_0),$ for all $(x,y)\in V_j.$ We define $\sigma _j:V_j\to Z^I$ an $(f,g)$-motion planner, given as
$$\sigma _j(x,y)(t):=\begin{cases}
\mathrm{pr}_1(F_j(x,y,2t)), & 0\leq t\leq \frac{1}{2} \\
\mathrm{pr}_2(F_j(x,y,2-2t)), & \frac{1}{2}\leq t\leq 1
\end{cases}$$
This proves that $\mathrm{TC}_H(f,g)\leq \mathrm{cat}(f\times g)$ and we conclude.
\end{proof}

\begin{corollary}
Let $X \xrightarrow{f} Z \xleftarrow{g} Y$ be a cospan with $Z$ a path-connected space. Then $\mathrm{TC}_H(f,g)=0$ if and only if $f$ and $g$ are nullhomotopic.
\end{corollary}
\begin{proof}
If $\mathrm{TC}_H(f,g)=0$, then
$$
\max\{\mathrm{cat}(f),\mathrm{cat}(g)\} \leq \mathrm{TC}_H(f,g)=0,
$$
so $\mathrm{cat}(f)=\mathrm{cat}(g)=0$, and both $f$ and $g$ are nullhomotopic.

Conversely, if $f$ and $g$ are nullhomotopic, then $f\times g$ is also nullhomotopic. Hence 
$\mathrm{TC}_H(f,g)\leq \mathrm{cat}(f\times g)=0,
$ and the claim follows.
\end{proof}

\begin{corollary}
Let $X \xrightarrow{f} Z \xleftarrow{g} Y$ be a cospan with $Z$ a path-connected space. If $f$ is nullhomotopic, then $\mathrm{TC}_H(f,g)=\mathrm{cat}(g).$ Similarly, if $g$ is nullhomotopic, then $\mathrm{TC}_H(f,g)=\mathrm{cat}(f).$
\end{corollary}

\begin{proof}
Just observe that, if $f$ is nullhomotopic, then $\mathrm{cat}(f)=0$ and $\mbox{cat}(f\times g)\leq \mbox{cat}(g)$. Applying Proposition \ref{gen-bounds} we obtain the result. Similarly when $g$ is nullhomotopic.    
\end{proof}
The next result establishes the behavior of $\mathrm{TC}_H(f,g)$  under pre- and post-composition.

\begin{proposition}\label{post-pre}
Let $X \xrightarrow{f} Z \xleftarrow{g} Y$ be a cospan.
\begin{enumerate}
\item If $u:X' \rightarrow X$ and $v:Y' \rightarrow Y$ are maps, then
$\mathrm{TC}_H(f\circ u,g\circ v)\leq \mathrm{TC}_H(f,g).$
\item If $w:Z\to Z'$ is a map, then $\mathrm{TC}_H(w\circ f,w\circ g)\leq \mathrm{TC}_H(f,g).$
\end{enumerate}

\end{proposition}
\begin{proof}
For (1), taking into account that $(f\times g)\circ (u\times v)=(f\circ u)\times (g\circ v)$ we have the following consecutive pullbacks
$$\xymatrix{
{E_{f\circ u,g\circ v}} \ar[d]_{\theta _{f\circ u,g\circ v}} \ar[r]  & {E_{f,g}} \ar[d]_{\theta _{f,g}} \ar[r]  & {Z^I} \ar[d]^{\pi _Z} \\
{X'\times Y'} \ar[r]_{u\times v}  & {X\times Y} \ar[r]_{f\times g}  & {Z\times Z}
}$$
Therefore $\mathrm{TC}_H(f\circ u,g\circ v)=\mathrm{secat}(\theta _{f\circ u,g\circ v})\leq \mathrm{secat}(\theta _{f,g})=\mathrm{TC}_H(f,g).$

For (2) just observe the existence of a canonical map $E_{f,g}\to E_{w\circ f,w\circ g}$ making commutative the following triangle
$$\xymatrix{
{E_{f,g}} \ar[rr] \ar[dr]_{\theta _{f,g}} & & {E_{w\circ f,w\circ g}} \ar[dl]^{\theta _{w\circ f,w\circ g}} \\ & {X\times Y} &
}$$ Therefore, by Proposition \ref{triangle}, $\mathrm{TC}_H(w\circ f,w\circ g)=\mathrm{secat}(\theta _{w\circ f,w\circ g})\leq \mathrm{secat}(\theta _{f,g})=\mathrm{TC}_H( f,g).$
\end{proof}

\begin{corollary}
Let $X \xrightarrow{f} Z \xleftarrow{g} Y$ be a cospan. Given the subspaces $A\subseteq X$ and $B\subseteq Y$ we have
$$\mathrm{TC}_H(f_{\vert A} ,g_{\vert B})\leq \mathrm{TC}_H(f,g)$$
\end{corollary}
\begin{proof}
Just apply Proposition \ref{post-pre} (1) when $u:A\hookrightarrow X$ and $v:B\hookrightarrow Y$ are the respective inclusions.
\end{proof}

\begin{corollary}
Let $X \xrightarrow{f} Z \xleftarrow{g} Y$ be a cospan. If $f$ and $g$ have a right homotopy inverse, then:
$$\mathrm{TC}_H(f,g)=\mathrm{TC}(Z).$$
In particular, 
$\mathrm{TC}_H(\mathrm{pr}_1,\mathrm{pr}_2)=\mathrm{TC}(X),$ where $\mathrm{pr}_1,  \mathrm{pr}_2:X\times X\to X$ denote the projection maps.
\end{corollary}

\begin{proof}
If $f'$ and $g'$ denote the respective right homotopy inverses of $f$ and $g$, then applying Proposition \ref{gen-bounds} and Proposition \ref{post-pre} (1) we obtain
$$\mathrm{TC}(Z)=\mathrm{TC}_H(id_Z,id_Z)=\mathrm{TC}_H(f\circ f',g\circ g')\leq \mathrm{TC}_H(f,g)\leq \mathrm{TC}(Z).$$
\end{proof}

To conclude our study, we consider a cospan of maps
$X \xrightarrow{f} Z \xleftarrow{g} Y,$
where $Z$ is an $H$-group. We denote by 
$\mu : Z \times Z \longrightarrow Z$
the multiplication, and by 
$\nu : Z \longrightarrow Z$
the homotopy inverse. We denote by 
$e$ the unit element, chosen as basepoint.  Associated with this structure we define the \emph{difference map}, $\delta :Z\times Z\to Z$, as $\delta (x,y):=\mu (x,\nu(y))$ and its composition with $f\times g$
$$\delta _{f,g}:X\times Y\to Z, \qquad 
\delta_{f,g}(x,y)=\mu\big(f(x),\nu(g(y))\big).$$
It is well known that $\mathrm{TC}(Z)=\mathrm{cat}(Z)$ and therefore, by Proposition \ref{gen-bounds}(1), we have $\mathrm{TC}_H(f,g)\leq \mathrm{cat}(Z).$ The next theorem gives a more accurate estimation of $\mathrm{TC}_H(f,g):$

\begin{theorem}\label{H-groups}
Let $X \xrightarrow{f} Z \xleftarrow{g} Y$ be a cospan of maps, where $Z$ is a path-connected $H$-group. Then
$$\mathrm{TC}_H(f,g)=\mathrm{cat}(\delta _{f,g}).$$

\end{theorem}

\begin{proof}
Suppose that $\mathrm{TC}_H(f,g)=n$ and take $\{U_i\}_{i=0}^n$ an open cover of $X\times Y$ with $(f,g)$-motion planners
$$\sigma _i:U_i\to Z^I$$
By the $H$-group structure of $Z$ we take a homotopy
$F:Z\times I\to Z$ such that $F(z,0)=\mu (z,\nu(z))$ and $F(z,1)=e$, for all $z\in Z.$ Then we define
$G_i:U_i\times I\to Z$ by
$$G_i(x,y,t):=\begin{cases}
\mu (\sigma _i(x,y)(2t),\nu (g(y))), & 0\leq t\leq \frac{1}{2} \\
F(g(y),2t-1), & \frac{1}{2}\leq t\leq 1
\end{cases}$$
This construction defines a nullhomotopy $\delta _{f,g}|U_i\simeq C_e$. Therefore, $\mathrm{cat}(\delta _{f,g})\leq n.$

Conversely, suppose $\mathrm{cat}(\delta _{f,g})=n$ and take $\{V_j\}_{j=0}^n$ an open cover of $X\times Y$ with nullhomotopies
$$F_j:V_j\times I\to Z$$ \noindent with $F_j(x,y,0)=\delta _{f,g}(x,y)=\mu (f(x),\nu (g(y)))$ and $F_j(x,y,1)=e,$ for all $(x,y)\in V_j.$ 
By a combination of homotopy associativity with the homotopy inverse property we can take a homotopy
$$L:Z\times Z\times I\to Z$$
\noindent such that $L(z,z',0)=\mu (\mu (z,\nu (z')),z')$ and $L(z,z',1)=z,$ for all $(z,z')\in Z\times Z.$ Moreover, we can also take, from the homotopy unit property, another homotopy
$$M:Z\times I\to Z$$ \noindent satisfying $M(z,0)=\mu (e,z)$ and $M(z,1)=z,$ for all $z\in Z.$ Then we define a $(f,g)$-motion planning $\sigma_j:V_j\to Z^I$ by the following formula
$$\sigma _j(x,y)(t):=\begin{cases}
L(f(x),g(y),1-3t), & 0\leq t\leq \frac{1}{3} \\

\mu (F_j(x,y,3t-1),g(y)), & \frac{1}{3}\leq t\leq \frac{2}{3} \\

M(g(y),3t-2), &  \frac{2}{3}\leq t\leq 1
\end{cases}$$
Therefore $\mathrm{TC}_H(f,g)\leq n$ and we conclude the proof.
\end{proof}

\begin{remark}
The same conclusion remains valid if $Z$ is assumed to be only a path-connected $H$-space endowed with a \emph{difference map} $\delta$. By this we mean a continuous map $\delta : Z \times Z \to Z$ such that
\[
\mu \circ (\mathrm{pr}_1, \delta) \simeq \mathrm{pr}_2
\qquad \text{and} \qquad
\delta \circ \Delta_Z \simeq C_e ,
\]
where $\mu$ is the $H$-space multiplication, $\Delta_Z : Z \to Z \times Z$ denotes the diagonal map, and $C_e : Z \to Z$ is the constant map at the unit element $e$.
\end{remark}

\section*{Discussion and Outlook}

This work develops a bivariate version of topological complexity associated with a pair 
of maps $X \xrightarrow{f} Z \xleftarrow{g} Y$. The resulting invariant extends both 
Farber’s classical topological complexity and Pavešić’s map-based formulation, and provides 
a systematic way to measure the cost of synchronizing two systems through a common target space. 
It fits naturally into the framework of sectional category and strict sectional number, and the 
constructions and proofs are developed consistently in this setting rather than through filtered versions.

The basic structure of the invariant has been established: symmetry in $f$ and $g$, subadditivity, behavior with respect to pre- and post-compositions, and comparison with the 
homotopy-invariant version $\mathrm{TC}_{H}(f,g)$. The collaboration principle shows that the presence 
of a fibration can substantially reduce the complexity of coordination. In addition, a cohomological 
lower bound has been obtained from the image of the zero-divisor ideal, extending the classical method 
used for univariate topological complexity.

The strict and homotopy-invariant formulations should be regarded as complementary. The strict version 
detects rigid synchronization requirements and depends on the specific representatives of the maps, 
while the homotopy-invariant version provides a homotopy-stable lower bound and remains finite under 
assumptions where the strict invariant does not. Their comparison clarifies the role of fibrations and 
homotopy invariance in the bivariate setting.

Beyond the results established here, several directions remain open. One natural problem is the refinement of cohomological estimates, possibly involving mixed algebraic structures arising from the interaction of $H^*(X)$, $H^*(Y)$, and $H^*(Z)$, as well as the development of a genuinely bivariate notion of Lusternik--Schnirelmann category compatible with the present framework.
From a different perspective, the bivariate approach suggests extensions to other settings in which sectional methods remain meaningful. In particular, simplicial variants require additional compatibility conditions to preserve the interpretation of synchronization and sectional obstructions, and cannot be obtained by a purely formal translation of the continuous theory. Preliminary work in this direction indicates the appearance of phenomena that are absent in the homotopy-invariant case.
Some of these developments are currently under investigation and will be reported elsewhere.

\end{document}